\documentclass[12pt]{article}

\usepackage{amsfonts,amsmath,amsthm,amssymb,amscd}
\usepackage[dvips]{graphicx}

\newtheorem{theorem}{Theorem}
\newtheorem{lemma}[theorem]{Lemma}
\newtheorem{proposition}[theorem]{Proposition}
\newtheorem{corollary}[theorem]{Corollary}

\numberwithin{equation}{section}
\numberwithin{theorem}{section}

\newcommand{\trace}{\qopname\relax o{tr}}

\newcommand{\diag}{\qopname\relax o{diag}}

\newcommand{\bbP}{\mathbb{P}}
\newcommand{\C}{\mathbb{C}}
\newcommand{\Z}{\mathbb{Z}}

\newcommand{\T}{\mathbb{T}}

\newcommand{\Res}{\textrm{Res}}

\begin{document}

\begin{center}
{\bf 
{\Large 
Isomonodromic deformation with \\
an irregular singularity 
and the $\Theta$ function
}
}
\end{center}

\begin{center}
 
Kazuhide Matsuda
\end{center}

\begin{center}
Department of Mathmatics, 
Graduate School 
of 
Information Science and Technology, \\
Osaka University,\\
1-1 Machikaneyama-machi, 
Toyonaka, 560-0043. 
\end{center}

\vskip 1mm
\par\noindent
{\bf Abstract:} \,
We will study 
a monodromy preserving deformation 
with 
an irregular singular point 
and 
determine the $\tau$ function 
of 
the monodromy preserving deformation 
by 
the elliptic $\theta$ function
moving 
the argument $z$ 
and 
the period $\Omega$. 

\section{Introduction}
In this paper, 
we construct 
a monodromy preserving deformation 
whose 
$\tau$ function 
is represented 
by the $\theta$ function 
moving 
both 
the argument $z$ 
and 
the period $\Omega$. 
The monodromy preserving deformation 
has an irregular singularity of 
the Poincar\'e rank one. 
\par
In \cite{MM1}, 
Miwa-Jimbo-Ueno 
extended 
the work of Schlesinger in \cite{Sch}
and 
established a general theory 
of monodromy preserving deformation 
for a first order matrix system 
of 
ordinary linear differential equations, 
\begin{equation}
\frac{dY}{dx}
=
A(x)Y,
\quad
A(x)
=
\sum_{\nu=1}^{n}
\sum_{k=0}^{r_{\nu}}
\frac{A_{\nu,-k}}{(x-a_{\nu})^{k+1}}
-
\sum_{k=1}^{r_{\infty}}
A_{\infty, -k} x^{k-1},
\end{equation}
having regular or irregular singularities 
of arbitrary rank.
\par
The monodromy data to be preserved is 
\newline
\newline
(i) Stokes multipliers $S_j^{(\nu)} \,\,
(j=1, \ldots , 2r_{\nu}),$
\newline
\newline
(ii) connection matrices $C^{(\nu)},$
\newline
\newline
(iii) ``exponents of formal monodromy" $T_0^{(\nu)}.$
\newline
\newline
Miwa-Jimbo 
found 
a deformation equation 
as a necessary and sufficient condition 
for the monodromy data 
to be invariant 
for deformation parameters 
and 
defined 
the $\tau$ function 
for the deformation equation.
\par
We will explain the relationship 
between 
the $\tau$ function 
and 
the $\theta$ function. 
In \cite{MM2}, 
Miwa-Jimbo 
expressed 
the $\tau$ function 
with 
the $\theta$ function 
by 
moving 
the argument $z$ 
of the $\theta$ function. 
The monodromy preserving deformation 
has 
irregular singularities. 
The Riemann surface 
is 
a ramified covering 
of 
$\mathbb{CP}^1$ 
with 
the 
covering degree $m$ 
and 
its genus is $g$. 
Especially 
the case of genus one is written 
in \cite{MM0}, 
which we will describe in Appendix.
\par
In \cite{KK}, 
Kitaev-Korotkin 
calculated the $\tau$ function 
with 
the $\theta$ function 
by 
moving the period $\Omega$ 
of the $\theta$ function. 
They consider monodromy preserving deformation of
a second order Fuchsian equation with 
$2g+2$ regular singular points.
In this case the $\tau$ function is represented 
by  the  theta  null value of 
a hyperelliptic curve with genus $g$. 
In case $g=1$, Kitaev-Korotkin's monodromy preserving deformation 
is equivalent to Picard's solution of the sixth Painlev\'e equation 
by the B\"acklund transformation. 
The aim of this paper 
is 
to unify Miwa-Jimbo's result 
and Kitaev-Korotkin's result 
in the elliptic case.
\par
It is well known that 
the elliptic $\theta$ function 
satisfies the heat equation, 
whose solution space is 
infinite dimensional. 
The equation (5.4) of 
isomonodromic deformation is 
solved by the $\tau$ function. 
The elliptic $\theta$ function 
also satisfies 
the nonlinear equation (5.4), 
whose solution space is finite dimensional.
\par
We will explain 
this paper in detail. 
In section 2, 
we show some definitions 
and facts of elliptic functions, 
which is necessary 
to construct the $\tau$-function. 
\par
In section 3, 
we will study 
the following ordinary differential equation:
\begin{equation}
\label{equ:diff1}
\frac{dY}{dx}
= 
\left( 
\frac{B_{-1}}{(x-a)^2}
+
\frac{B_0}{x-a}
+
\sum_{\nu=1}^3
\frac{A_{\nu}}{x-e_{\nu}}
\right)
Y(x),
\end{equation}
and define the $\tau$-function.
\par
In section 4, 
we construct 
a fundamental solution of (\ref{equ:diff1}) 
and 
calculate the following monodromy data: 
\newline
\newline
(i) Stokes multipliers around $x=a$ \,\,
$S^a_1 = S^a_2 =1,$
\newline
\newline
(ii) connection matrices around $x=e_{\nu}$ \,\,
$C^{(\nu)} \,(\nu=1,2,3,\infty),$
\newline
\newline
(iii) the exponents of formal monodromy 
$
T^a_0=0, T^{(\nu)}_0=\diag(-\frac14, \frac14) 
\, (\nu=1,2,3,\infty),
$
\newline
\newline
where $e_{\infty}$ is $\infty$.
\par 
The behavior of the fundamental solution 
near the singular points 
is as follows 
\begin{eqnarray}
Y(x) 
&=&
(
1
+
O(x-a)
)
\exp
T^{(a)}(x), \\
Y(x)
&=&
G^{(\nu)}
(
1
+
O(x-e_{\nu})
)
\exp T^{(\nu)}(x) \,\,
(\nu=1,2,3), \\
Y(x)
&=&
G^{(\infty)}
(
1
+
O(\frac{1}{x})
)
\exp
T^{(\infty)}(x)
\end{eqnarray}
where
\begin{eqnarray}
T^{(a)}(x)
&=&
\left(
\begin{array}{cc}
\frac{\wp^{\prime}(\alpha)t}{2} &     \\
                 & -\frac{\wp^{\prime}(\alpha)t}{2}  \\
\end{array}
\right)
\frac{(x-a)^{-1}}{-1}
+
T_0^{(a)}
\log(x-a),  \\
T^{(\nu)}(x)
&=&
T^{(\nu)}_0
\log(x-e_{\nu})
\,\, (\nu=1,2,3),  \\
T^{(\infty)}(x)
&=& 
T_0^{(\infty)}
\log
(
\frac{1}{x}
).
\end{eqnarray}
\par
In section 5, 
we concretely calculate 
the monodromy preserving deformation 
from the fundamental solution 
and compute the coefficients, 
$B_{-1}, B_0, A_{\nu} \,\,(\nu=1,2,3).$ 
\par
In section 6, we calculate the $\tau$-function. 
Section 6 consists of three subsections. 
Subsection 6.1 is devoted to $H_t$, 
the Hamiltonian on the deformation parameter $t$.
Subsection 6.2, 6.3 is about 
$H_{\nu} \,(\nu=1,2,3)$, 
the Hamiltonians on the deformation parameters 
$e_{\nu} \, (\nu=1,2,3).$ 
In subsection 6.2, 
we show some facts about 
elliptic functions 
in order to calculate the Hamiltonians 
$H_{\nu} \,(\nu=1,2,3)$. 
In subsection 6.3, 
we calculate $H_{\nu} \,(\nu=1,2,3)$ 
and 
the $\tau$-function.
Our main theorem is as follows: 
\begin{theorem}
For the monodromy preserving deformation 
(\ref{equ:diff1}),
the $\tau$ function is 
\begin{eqnarray*}
\tau(e_1, e_2, e_3 , t)
&=& 
\theta[p,q]
\left(
\frac{t}{\omega_1}
;\Omega
\right)
\omega_1^{-\frac12}
\prod_{1 \leq \nu <\mu \leq 3} 
(e_{\nu}-e_{\mu})^{-\frac18}   
\exp(\frac{\eta_1 t^2}{2 \omega_1})
\exp(\frac{t^2}{4} f(e_1, e_2, e_3)),
\end{eqnarray*}
where 
\begin{eqnarray*}
f(e_1, e_2, e_3) 
&=& 
-a 
+ 
\frac13 \sum_{\nu=1}^3 e_{\nu}  
+ 
\frac12 \prod_{\nu=1}^3(a-e_{\nu})
\sum_{\mu=1}^3 \frac{1}{(a-e_{\mu})^2}.
\end{eqnarray*}
\end{theorem}
\par\noindent
{\bf Acknowledgments.} 
The author thanks 
Professor Yousuke Ohyama 
for the careful guidance. 
This work is partially done 
by the author 
in the Isaac Newton Institute 
for Mathematical sciences, 
Cambridge University. 
The author is grateful 
to Professor D. Korotkin 
for a fruitful discussion.

\section{Elliptic Function Theory}
We will review the elliptic function theory 
to fix our notation. 
See \cite{Mum}. 
We will consider the elliptic curve $E$ 
defined by the equation
$$
y^2 = 4(x-e_1)(x-e_2)(x-e_3)
$$
with arbitrary constants 
$e_i \in \C \,\,(e_i \neq e_j),$ 
where we do not assume that 
$ e_1+e_2+e_3=0$. 
We set the new coordinate:
$$
\tilde{x}
=x
-
\frac{1}{3} \sum_{i=1}^3 e_i, \,\, 
\tilde{e_j}
=
e_j - \frac13 \sum_{i=1}^3 e_i 
\,\,\,
\textrm{for} \,\,j=1,2,3 
$$
and define the Abel map in the following way:
\begin{equation}
u
:= 
\int_{\infty}^x 
\frac{dx}{\sqrt{4(x-e_1)(x-e_2)(x-e_3)}}
=
\int_{\infty}^{\tilde{x}} 
\frac
{
d\tilde{x}
}
{
\sqrt
{
4
(\tilde{x}-\tilde{e}_1)
(\tilde{x}-\tilde{e}_2)
(\tilde{x}-\tilde{e}_3)
}
}.
\end{equation}
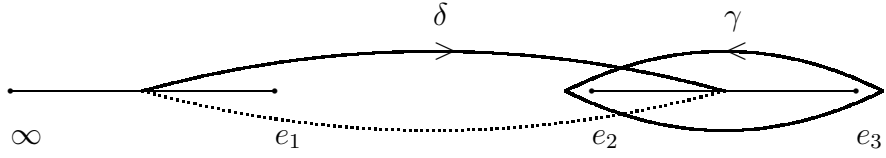
\begin{figure}
\begin{center}
\begin{picture}(400,90)
\put(40,30){\circle*{2}}
\put(40,10){$\infty$ }
\put(140,30){\circle*{2}}
\put(140,10){$e_1$}
\thinlines
\put(40,30){\line(1,0){100}}

\thicklines
\qbezier(90,30)(200,60)(310,30)

\put(200,42){$>$}
\put(200,55){$\delta$}

\thicklines
\qbezier[100](90,30)(200,0)(310,30)

\put(260,30){\circle*{2}}
\put(260,10){$e_2$}
\put(360,30){\circle*{2}}
\put(360,10){$e_3$}
\thinlines
\put(260,30){\line(1,0){100}}

\thicklines
\qbezier(250,30)(310,60)(370,30)

\put(310,55){$\gamma$}

\put(310,42){$<$}

\thicklines
\qbezier(250,30)(310,0)(370,30)

\end{picture}
\end{center}
\caption
{
Branch cuts 
and 
canonical basis 
of cycles on the elliptic curve, 
$E$. 
Continuous 
and 
dashed paths 
lie 
on the first and second 
sheet of $E$, respectively.
}
\label{rs}
\end{figure}

Let us 
define periods $\omega_1, \omega_2$ by
\begin{equation}
\omega_1
=
\int_{\gamma} \frac{dx}{y}, 
\,\, 
\omega_2=\int_{\delta} \frac{dx}{y}.
\end{equation}
From the bilinear relation 
of Riemann, we get 
\begin{equation}
\Omega
:=
\frac{\omega_2}{\omega_1}, 
\,\, 
\Im (\Omega) >0.
\end{equation}
$ 
\tilde{x} 
= 
\tilde{x}(u)
$, 
the inverse function of the Abel map, 
can be expressed 
with the Weierstrass $\wp$-function:
\begin{equation}
x
-
\frac{1}{3} 
\sum_{i=1}^3 e_i
=
\tilde{x}(u) 
= 
\wp(u; \omega_1, \omega_2).
\end{equation} 
$\wp(u)$ satisfies 
the following differential equations:
\begin{eqnarray}
\wp^{\prime}(u)^2
&=&
4
(\wp(u)-\tilde{e_1})
(\wp(u)-\tilde{e_2})
(\wp(u)-\tilde{e_3})  \\
&=&
4 
\wp(u)^3 
- 
g_2 
\wp(u) 
- 
g_3
\end{eqnarray}
We have
\begin{equation}
\label{eqn:sig}
\sigma
(u; 
\omega_1, 
\omega_2)
= 
\exp 
\left( 
\frac{\eta_1}{2\omega_1}u^2 
\right) 
\frac{\omega_1}{\theta_{11}^{\prime}} 
\theta_{11} 
\left( 
\frac{u}{2\omega_1}; \Omega 
\right).
\end{equation}
We set 
the $\sigma$ function with characteristic 
$p, q \in \C$ in the following way:
\begin{equation}
\sigma [p,q] (u) 
:= 
\exp 
\left( 
\frac{\eta_1}{2 \omega_1} u^2 
\right)
\frac{\omega_1}{\theta_{11}^{\prime}}
\theta [p,q] 
\left( 
\frac{u}{\omega_1}; \Omega 
\right).
\end{equation}
$\sigma [p,q] (u)$ 
has 
the following periodicity properties:

\begin{equation}
\label{eqn:per1}
\sigma[p,q](u+\omega_1) 
= 
\exp (2\pi i p) 
\exp \eta_1 
\left(
u
+
\frac{\omega_1}{2} 
\right)
\sigma[p,q] (u), 
\end{equation}
\begin{equation}
\label{eqn:per2}
\sigma[p,q](u+\omega_2) 
= 
\exp (-2 \pi iq) 
\exp \eta_2 
\left(
u
+
\frac{\omega_2}{2} 
\right) 
\sigma[p.q](u).
\end{equation}

\section{The Schlesinger System}
We will 
study the following ordinary differential equation:
\begin{equation}
\label{equ:diff2}
\frac{dY}{dx}
= 
\left( 
\frac{B_{-1}}{(x-a)^2}
+
\frac{B_0}{x-a}
+
\sum_{{\nu}=1}^3
\frac{A_{\nu}}{x-e_{\nu}}
\right)
Y(x),
\end{equation}
where $A_1, A_2, A_3, B_{-1}, B_0 \in sl(2, \C)$ 
are independent of $x$. 
\par
The monodromy data of \eqref{equ:diff2} 
is as follows:
\newline
\newline
(i) Stokes multipliers around $x=a$ \,\, 
$S^a_1 = S^a_2 ;$
\newline
\newline
(ii) connection matrices around $\infty, e_{\nu}$ \,\,
$C^{\infty}, C^{(\nu)} \,(\nu=1,2,3);$
\newline
\newline
(iii) the exponents of formal monodromy 
$T^a_0, T^{(\infty)}_0, T^{(\nu)}_0 \, (\nu=1,2,3).$
\newline
\newline
In the next section, 
we will get 
convergent series around 
$x = a, \infty, e_{\nu} (\nu=1,2,3)$ 
in the following way:
\begin{eqnarray}
& &
Y(x) 
=
(
1
+
O(x-a)
)
\exp
T^{(a)}(x)
=\hat{Y}^{(a)}(x) 
\exp T^{(a)}(x), \\
& &
Y(x)
=
G^{(\infty)}
(
1
+
O(\frac{1}{x})
)
\exp
T^{(\infty)}(x)
=G^{(\infty)}
\hat{Y}^{(\infty)}(x) 
\exp 
T^{(\infty)}(x), \\
& &
Y(x)
=
G^{(\nu)}
(
1
+
O(x-e_{\nu})
)
\exp T^{(\nu)}(x) 
=
G^{(\nu)}
\hat{Y}^{(\nu)}(x)
\exp
T^{(\nu)}_0(x),  \notag \\
& & 
\hspace{90mm}(\nu=1,2,3).
\end{eqnarray}
where
\begin{eqnarray}
T^{(a)}(x)
&=&
\left(
\begin{array}{cc}
\frac{\wp^{\prime}(\alpha)t}{2} &     \\
                  & -\frac{\wp^{\prime}(\alpha)t}{2}  \\
\end{array}
\right)
\frac{(x-a)^{-1}}{-1}
+
T_0^{(a)}
\log(x-a),  \\
T^{(\infty)}(x)
&=& 
T_0^{\infty}
\log
(
\frac{1}{x}
) \\
T^{(\nu)}(x)
&=&
T^{(\nu)}_0
\log(x-e_{\nu})
\,\, (\nu=1,2,3).
\end{eqnarray}

\par
For the deformation parameters 
$t, e_1, e_2, e_3,$ 
we set the following closed 1-form
\begin{eqnarray}
\omega 
&=& 
\omega_a 
+
\omega_{e_1}
+
\omega_{e_2}
+
\omega_{e_3} \\
&=& 
H_t dt 
+ 
H_1 d e_1 
+ 
H_2 d e_2 
+ 
H_3 d e_3, 
\end{eqnarray}
where
\begin{eqnarray}
\omega_{a}
&=& 
- 
\Res_{x=a} \,\, 
\trace
\hat{Y}^{(a)}(x)^{-1} 
\frac{\partial \hat{Y}^{(a)}}{\partial x}
(x)
\, 
d T^{(a)}(x)  ,\\
\omega_{e_{\nu}}
&=& 
- 
\Res_{x=e_{\nu}} \,\, 
\trace
\hat{Y}^{(\nu)}(x)^{-1} 
\frac{\partial \hat{Y}^{(\nu)}}{\partial x}(x)
\, 
d T^{(\nu)}(x) \,\,(\nu=1,2,3),
\end{eqnarray}
and
$d$ 
is the exterior differentiation 
with respect to 
the deformation parameters $t, e_1, e_2, e_3$. 
Especially, we can write
\begin{equation}
\omega_{e_{\nu}} 
= 
\left[ 
\Res_{x=e_{\nu}} 
\frac12 
\trace 
\left( 
\frac{dY}{dx}Y^{-1} 
\right)^2 
\right]
d e_{\nu}.
\end{equation}
From the closed 1-form $\omega$, 
we define the $\tau$-function in the following way:
\begin{equation}
\omega := d \log \tau (e_1, e_2, e_3,t).
\end{equation}

\section{The Riemann-Hilbert Problem for Special Parameters}
\quad 
In this section, 
we will concretely construct 
a $2\times2$ matrix valued function $Y(x)$, 
whose monodromy data, 
(i) Stokes multipliers, 
(ii) connection matrices, 
(iii) exponents of formal monodromies, 
is independent of 
the deformation parameters $t, e_1, e_2, e_3$.

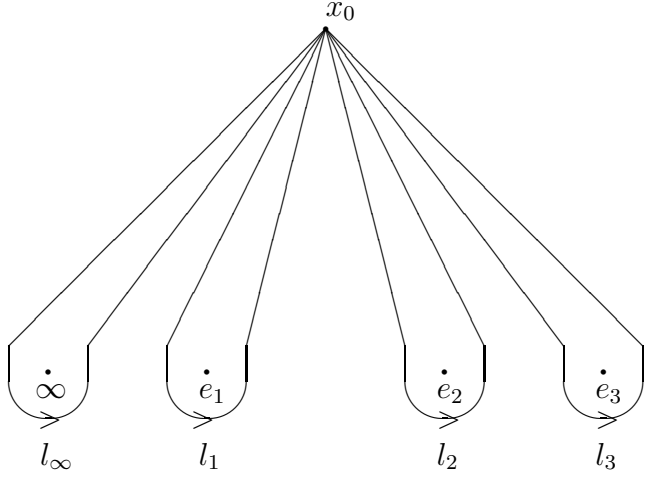
\begin{figure}
\begin{center}
\begin{picture}(400,150)

\put(70,70){\line(1,1){120}}
\put(100,70){\line(3,4){90}}

\put(85,70){\oval(30,55)[b]}
\put(85,60){\circle*{2}}
\put(80,50){$\infty$}
\put(81,39){$>$}
\put(82,25){$l_{\infty}$}

\put(130,70){\line(1,2){60}}
\put(160,70){\line(1,4){30}}

\put(145,70){\oval(30,55)[b]}
\put(145,60){\circle*{2}}
\put(142,50){$e_1$}
\put(139,39){$>$}
\put(142,25){$l_1$}

\put(190,190){\circle*{2}}
\put(190,195){$x_0$}

\put(220,70){\line(-1,4){30}}
\put(250,70){\line(-1,2){60}}

\put(235,70){\oval(30,55)[b]}
\put(235,60){\circle*{2}}
\put(232,50){$e_2$}
\put(232,39){$>$}
\put(232,25){$l_2$}

\put(280,70){\line(-3,4){90}}
\put(310,70){\line(-1,1){120}}

\put(295,70){\oval(30,55)[b]}
\put(295,60){\circle*{2}}
\put(292,50){$e_3$}
\put(292,39){$>$}
\put(292,25){$l_3$}

\end{picture}
\end{center}
\caption
{
Generators of 
$
\pi_1 
\left(
\C \bbP^1
\setminus
\{
\infty,
e_1,
e_2,
e_3
\},
x_0
\right)
$
}
\end{figure}

\par
We denote the Abel map $u=u(x)$ by $u=u(P)$, 
because we regard it as a function of 
the Riemann surface $E$ 
defined by
\begin{equation}
y^2= 4(x-e_1)(x-e_2)(x-e_3).
\end{equation}
We set the involution of $E$ in the following way:
\begin{equation}
*:
(x,y)
\longrightarrow
(x,-y).
\end{equation}
\par
We define the $2 \times 2$ matrix valued function 
$\Phi(P)$ 
with the $\sigma$ function 
and 
the $\zeta$ function in the following way:
\begin{equation}
\Phi (P) = 
\left(
\begin{array}{cc}
\varphi(P) & \varphi(P^{*}) \\
\psi(P) & \psi(P^{*})
\end{array}
\right),
\end{equation}
where
\begin{eqnarray*}
\varphi(P) 
&=& 
\sigma [p,q] (u+u_{\varphi}+t) 
\sigma(u-u_{\varphi})
\exp 
-\frac{t}{2} 
\{ 
\zeta(u-\alpha)
+
\zeta(u+\alpha) 
\}  \\
\psi(P) 
&=& 
\sigma [p,q] (u+u_{\psi}+t) 
\sigma(u-u_{\psi})
\exp 
-\frac{t}{2} 
\{ 
\zeta(u-\alpha)
+
\zeta(u+\alpha) 
\} \\
u_{\varphi} 
&=& 
u(P_{\varphi}), 
\, 
u_{\psi}
=u(P_{\psi}), 
\, 
\alpha
=
u(a),
\end{eqnarray*}
with arbitrary $P_{\varphi}, P_{\psi} \in E$.
\begin{proposition}
\label{prop:reg}
The function $\Phi(P)$ is holomorphic 
and 
invertible outside of 
the branch points $\infty, e_1, e_2, e_3$ 
and 
transforms as follows with respect to 
the tracing along the basic cycles
$\gamma, \delta$:
\begin{eqnarray*}
T_{\gamma} 
\big[ 
\Phi (P) 
\big] 
&=&
\Phi (P) 
\left(
\begin{array}{cc}
\exp \{ \pi i (2p+1) \} &        \\ 
           & \exp \{ - \pi i (2p+1) \} 
\end{array}
\right) 
\exp \{ \eta_1 (2u+\omega_1) \} \\
T_{\delta} 
\big[ 
\Phi (P) 
\big] 
&=&
\Phi (P) 
\left(
\begin{array}{cc}
\exp \{ - \pi i (2q+1) \} &         \\
              &  \exp \{ \pi i (2q+1) \}
\end{array}
\right)
\exp \{ \eta_2 (2u+\omega_2) \}
\end{eqnarray*} 
where by $T_{l}$ 
we denote 
the operator of analytic continuation 
along the contour $l$.
\end{proposition}
\begin{proof}
By using the periodicity 
(\ref{eqn:per1}) 
and 
(\ref{eqn:per2}), 
we get 
\begin{eqnarray*}
T_{\gamma} 
\left[ 
\Phi (P) 
\right] 
&=&
\Phi (P) 
\left(
\begin{array}{cc}
\exp \{ \pi i (2p+1) \} &        \\ 
                & \exp \{ - \pi i (2p+1) \} 
\end{array}
\right) 
\exp \{ \eta_1 (2u+\omega_1) \} \\
T_{\delta} 
\left[ 
\Phi (P) 
\right] 
&=&
\Phi (P) 
\left(
\begin{array}{cc}
\exp \{ - \pi i (2q+1) \} &        \\
             &  \exp \{ \pi i (2q+1) \}
\end{array}
\right)
\exp \{ \eta_2 (2u+\omega_2) \}.
\end{eqnarray*}
Then we obtain 
\begin{eqnarray*}
T_{\gamma}
\left[
\det \Phi (P)
\right]
&=&
\exp 2 \eta_1 (2 u + \omega_1)
\det \Phi (P)  \\
T_{\delta}
\left[
\det \Phi(P)
\right]
&=&
\exp 2 \eta_2 (2 u + \omega_2)
\det \Phi(P).
\end{eqnarray*}
By integrating 
$
\frac
{
d(\det \Phi(P)) 
}
{\det \Phi(P)}
$
along the fundamental polygon of $E$ 
and using the Legendre relation, we get
\begin{equation}
\frac{1}{2 \pi i}
\oint_{\partial \hat{E}}
\frac
{
d(\det \Phi(P))
}
{\det \Phi(P)}
=
\frac{1}{2 \pi i}
4(\eta_1 \omega_2 - \eta_2 \omega_1)
=
4.
\end{equation}
Therefore, four zeros of $\det \Phi(P)$ on $E$ 
are branch points, 
$\infty, e_1, e_2, e_3$. 
We proved the proposition.
\end{proof}
We will normalize $\Phi(P)$ 
at $x=a.$ 
We can easily show 
the local behavior 
of 
the Abel map.

\begin{lemma}
\label{lem:inv}
\begin{equation}
u-\alpha
=
\frac{1}{\wp^{\prime}(\alpha)}
(x-a)
-
\frac{\wp^{\prime \prime}(\alpha)}
{
2 \wp^{\prime \prime \prime}(\alpha)
}
(x-a)^2
+
\left(
\frac{\wp^{\prime \prime}(\alpha)^2}
{
2 \wp^{\prime}(\alpha)^5
}
-
\frac
{
\wp^{\prime \prime \prime}(\alpha)
}
{6 \wp^{\prime}(\alpha)^4}
\right)
(x-a)^3
+
\cdots.
\end{equation}
\end{lemma}

From Lemma \ref{lem:inv}, 
$\Phi(x)$ can be developed near $x=a$ 
in the following way:
\begin{equation*}
\Phi(x)
=
\left(
G^{(a)}
+
O(x-a)
\right)
\exp
\left(
\begin{array}{cc}
- 
\frac{\wp^{\prime}(\alpha)t}{2(x-a)} &     \\
                    & \frac{\wp^{\prime}(\alpha)t}{2(x-a)}  \\
\end{array}
\right),
\end{equation*}
where $G^{(a)}$ is 
a $2 \times 2$ matrix with the matrix elements:
\begin{eqnarray*}
(G^{(a)})_{11}
&=&
\sigma[p,q](\alpha+u_{\varphi}+t)
\sigma(\alpha-u_{\varphi})
\exp 
\Big\{
- \frac{t}{2}
\left(
\frac
{
\wp^{\prime \prime}(\alpha)
}
{
2 \wp^{\prime}(\alpha)
}
+
\zeta(2\alpha)
\right)
\Big\},  \\
(G^{(a)})_{21}
&=& 
\sigma[p,q](\alpha + u_{\psi} +t) 
\sigma(\alpha-u_{\psi}) 
\exp 
\Big\{
-
\frac{t}{2}
\left(
\frac
{
\wp^{\prime \prime}(\alpha)
}
{
2 \wp^{\prime}(\alpha)
} 
+ 
\zeta(2\alpha)
\right) 
\Big\}, \\
(G^{(a)})_{12}
&=&
\sigma[p,q](-\alpha +u_{\varphi} +t) 
\sigma(-\alpha -u_{\varphi})
\exp 
\Big\{
\frac{t}{2}
\left(
\frac
{
\wp^{\prime \prime}(\alpha)
}
{
2 \wp^{\prime}(\alpha)
} 
+
\zeta(2\alpha)
\right)
\Big\},   \\
(G^{(a)})_{22}
&=&
\sigma[p,q](-\alpha -u_{\psi} +t) 
\sigma(-\alpha-u_{\psi})
\exp 
\Big\{
\frac{t}{2}
\left(
\frac
{
\wp^{\prime \prime}(\alpha)
}
{
2 \wp^{\prime}(\alpha)
} 
+ 
\zeta(2\alpha)
\right)
\Big\}.
\end{eqnarray*}
From Proposition \ref{prop:reg},
\begin{equation*}
\det G^{(a)}
=
\det \Phi (a)
\neq
0.
\end{equation*}
We define a matrix valued function $Y(x)$ 
in the following way:
\begin{equation}
Y(P)
=
\frac
{
\sqrt{\det \Phi(a)}
}
{
\sqrt{\det \Phi(P)}
} 
(G^{(a)})^{-1} 
\Phi(P),
\end{equation}

By the normalization near $x=a$, 
we obtain the following lemma:

\begin{lemma}
\label{lem:ya}
\begin{eqnarray*}
& &
Y(x)
=
\left(
\left(
\begin{array}{cc}
1 & 0 \\
0 & 1   
\end{array}
\right)
+
O(x-a)
\right)
\exp 
T^{(a)}(x) \\
& &
T^{(a)}(x)
=
T^{(a)}_{-1}
\frac{(x-a)^{-1}}{-1}, \,\, 
T^{(a)}_{-1}
=
\left(
\begin{array}{cc}
\frac{\wp^{\prime}(\alpha)t}{2} &      \\
                 & \frac{-\wp^{\prime}(\alpha)t}{2}
\end{array}
\right).
\end{eqnarray*}
Especially, if 
$u_{\varphi}=\alpha, u_{\psi}=-\alpha,$ 
\begin{eqnarray*}
Y(x)
&=&
\left(
\left(
\begin{array}{cc}
1 & 0 \\
0 & 1
\end{array}
\right)
+
Y_1^{(a)}
(x-a)
+
\cdots
\right)
\exp T^{(a)}(x) \\
\left(Y_1^{(a)}\right)_{11}
&=&
\frac{1}{\wp^{\prime}(\alpha)}
\Big\{
\frac{\sigma[p,q]^{\prime}(t)}{\sigma[p,q](t)}
-
\frac{t}{2}
\left(
4\wp(\alpha)
-
\frac12
\left(
\frac
{
\wp^{\prime \prime}(\alpha)
}
{\wp^{\prime}(\alpha)}
\right)^2
\right)
\Big\} \\
\left(Y_1^{(a)}\right)_{21}
&=&
\frac
{
\sigma[p,q](2\alpha+t)
}
{
\sigma[p,q](t)
\sigma[p,q](2\alpha)
\wp^{\prime}(\alpha)
}
\exp
\Big\{
-t
\left(
\frac
{
\wp^{\prime \prime}(\alpha)
}
{
2\wp^{\prime}(\alpha)
}
+
\zeta(2\alpha)
\right)
\Big\} \\
\left(Y_1^{(a)}\right)_{12}
&=&
\frac
{
-
\sigma[p,q](-2\alpha+t)
}
{
\sigma[p,q](t)
\sigma[p,q](2\alpha)
\wp^{\prime}(\alpha)
}
\exp
\Big\{
t
\left(
\frac
{
\wp^{\prime \prime}(\alpha)
}
{
2\wp^{\prime}(\alpha)
}
+
\zeta(2\alpha)
\right)
\Big\} \\
\left(Y_1^{(a)}\right)_{22}
&=&
\frac{-1}{\wp^{\prime}(\alpha)}
\Big\{
\frac{\sigma[p,q]^{\prime}(t)}{\sigma[p,q](t)}
-
\frac{t}{2}
\left(
4\wp(\alpha)
-
\frac12
\left(
\frac
{
\wp^{\prime \prime}(\alpha)
}
{\wp^{\prime}(\alpha)}
\right)^2
\right)
\Big\}.
\end{eqnarray*}
\end{lemma}

$Y(x)$ has 
the following monodromy 
and 
Stokes matrices.
\begin{theorem}
\label{thm:mono}
For $\nu = \infty, 1, 2, 3,$ we set 
\begin{equation}
M_{\nu} 
=
\left(
\begin{array}{cc}
0                 & m_{\nu} \\
- m_{\nu}^{-1} & 0
\end{array}
\right).
\end{equation}
The matrix elements 
corresponding to the contour 
$
l_{\nu} \, 
(\nu=\infty, 1,2,3)
$ 
are 
given by the expressions
\begin{eqnarray*}
& & 
m_{\infty}  = -i, \,\, 
m_{1}     = i \exp (-2 \pi i p)  \\
& & m_{2} = -i \exp 2 \pi i (-p+q) \\
& & m_{3} = i \exp (2\pi i q).
\end{eqnarray*}
Stokes matrices are as follows
$$
S_1 
=
S_2
=
\left(
\begin{array}{cc}
1 & 0  \\
0 & 1 \\
\end{array}
\right).
$$
\end{theorem}

\begin{proof}
From Proposition \ref{prop:reg}, 
we note that
the branch points $\infty, e_1, e_2, e_3$ on $E$ are 
zero points of $\det \Phi(P)$ of the first order.
According to Proposition \ref{prop:reg}, we obtain
\begin{eqnarray*}
T_{\gamma}
\left[
Y(x)
\right]
&=&
Y(x) M_{3} M_{2}  \\
&=&
Y(x)
\left(
\begin{array}{cc}
\exp \{2 \pi i p \} &       \\
              & \exp \{- 2\pi i p \}
\end{array}
\right),
\end{eqnarray*} 
\begin{eqnarray*}
T_{\delta^{-1}}
\left[
Y(x)
\right]
&=&
Y(x) M_{2} M_{1}  \\
&=&
Y(x)
\left(
\begin{array}{cc}
\exp \{2 \pi i q\} &       \\
           & \exp \{-2 \pi i q\} \\
\end{array}
\right).
\end{eqnarray*}
\par
By the involution *, we get
$$
Y(x) M_{\infty} 
=
Y(x) 
\left(
\begin{array}{cc}
  & i \\
i &     
\end{array}
\right),
\,
\textrm{or}
\,
=
Y(x) 
\left(
\begin{array}{cc}
   & -i \\
-i &     
\end{array}
\right).
$$
We set 
\begin{eqnarray*}
M_{\infty}
=
\left(
\begin{array}{cc}
                  & m_{\infty}  \\
-m_{\infty}^{-1}  &
\end{array}
\right)
=
\left(
\begin{array}{cc}
    & -i  \\
-i  &
\end{array}
\right).
\end{eqnarray*}
\par
From the property of the monodromy group, we obtain
\begin{equation}
\label{eqn:mono}
M_{3} M_{2} M_{1} = M_{\infty}^{-1}.
\end{equation}
From 
\begin{eqnarray*}
M_{3} M_{2} 
&=& 
\diag 
(
\exp \{2 \pi i p\}, \exp \{ -2 \pi i p\}
), \\
M_{2} M_{1}
&=&
\diag
(
\exp \{2 \pi i q \}, \exp \{- 2 \pi i q \}
),
\end{eqnarray*}
we get 
\begin{eqnarray*}
M_{1}
&=&
\left(
\begin{array}{cc}
             & -m_{\infty} \exp \{-2 \pi i p \} \\
m_{\infty}^{-1} \exp \{2 \pi i p \}  &  
\end{array}
\right),  \\
M_{3}
&=&
\left(
\begin{array}{cc}
              &  -m_{\infty} \exp \{2 \pi i q \}  \\
m_{\infty}^{-1} \exp \{-2 \pi i q \}  &
\end{array}
\right).
\end{eqnarray*}
Therefore, we get
\begin{equation*}
M_{2}
=
\left(
\begin{array}{cc}
                &  m_{\infty} \exp \{2 \pi i (q-p) \}  \\
-m_{\infty}^{-1} \exp \{2 \pi i (p-q) \}  &   
\end{array}
\right),
\end{equation*}
from (\ref{eqn:mono}).
\par
From Lemma \ref{lem:ya}, the Stokes matrices are
$$
S_1=S_2=
\left(
\begin{array}{cc}
1 & 0 \\
0 & 1 
\end{array}
\right).
$$
\end{proof}
We can describe the monodromy data of $Y(x)$ 
in the following way.

\begin{corollary}
$Y(x)$ has the following monodromy data:
\newline
\newline
(i) Stokes multipliers $S^a_1=S^a_2=1,$
\newline
\newline
(ii) connection matrices 
$
C^{(\nu)}
=
\frac{1}{\sqrt{2 i m_{\nu}}}
\left(
\begin{array}{cc}
i & - m_{\nu} \\
i &   m_{\nu} 
\end{array}
\right)
\,\,
(\nu = \infty, 1, 2, 3),
$
\newline
\newline
(iii) formal monodromies 
$
T_0^{(\nu)}
=
\diag(-\frac14, \frac14)
\,\,
(\nu= \infty, 1, 2, 3).
$
\newline
\newline
Especially, 
the behaviors of $Y(x)$ 
near 
$e_{\nu} \, (\nu= \infty, 1, 2, 3)$ 
are 
\begin{eqnarray*}
& &
Y(x)
=
G^{(\nu)}
(
1
+
O(x_{\nu})
)
\exp T^{(\nu)} (x) 
C^{(\nu)} \\
& &
T^{(\nu)}(x)
=
T^{(\nu)}_0(x) \log x_{\nu} \\
& &
x_{\infty} = x^{-1}, \,
x_{\nu} = x-e_{\nu}
\end{eqnarray*}
\end{corollary}

\begin{proof}
(i) is clear. 
(ii) and (iii) can be obtained 
by diagonalizing the monodromy matrices 
$M_{\nu} \,\, (\nu=\infty, 1, 2, 3).$
\end{proof}

This corollary 
means 
that 
the monodromy data of $Y(x)$ 
is independent of 
the deformation parameters $t, e_1, e_2, e_3.$

\section{Monodromy Preserving Deformation}
In this section, 
we will prove that 
$Y(x)$ satisfies an ordinary differential 
equation of the following type:
\begin{equation*}
\frac{d Y}{d x}
=
\left(
\frac{B_{-1}}{(x-a)^2}
+
\frac{B_0}{x-a}
+
\sum_{\nu=1}^3
\frac{A_{\nu}}{x-e_{\nu}}
\right)
Y(x),
\end{equation*}
and 
concretely 
determine 
the coefficients 
$B_{-1}, B_0, A_{\nu} \,\, (\nu=1,2,3).$
\par
Firstly, we define 
symbols 
in order to 
describe the monodromy preserving deformation. 
For $\nu=1, 2, 3,$ 
we set 
\begin{eqnarray*}
D^{(\nu)}(u) 
&=&
\frac{2m_{\nu}}{m_{\infty}}
\varphi(u)\psi(u)
\left(
\frac{d}{du} \log \varphi(u) 
- 
\frac{d}{du} \log \psi(u)
\right), \\
\tilde{\omega}_1 
&=& 
\frac{\omega_1}{2}, \, 
\tilde{\omega}_2
=
\frac{\omega_1+\omega_2}{2}, \,
\tilde{\omega}_3
=
\frac{\omega_2}{2},   \\
\tilde{\eta}_1
&=&
\eta_1, \,
\tilde{\eta}_2
=
\eta_1
+
\eta_2, \,
\tilde{\eta}_3
=
\eta_2.
\end{eqnarray*}
We note 
that 
the monodromy matrices 
and 
the Stoke matrices 
do not have the parameters $P_{\varphi}, P_{\phi}$ 
in Theorem \ref{thm:mono}.
Then in this section, 
we set 
$
u_{\varphi}=\alpha, 
\, 
u_{\phi}=-\alpha
$ 
from the uniqueness of the Riemann-Hilbert problem.
\begin{theorem}
\label{thm:mpd}
$Y(x)$ satisfies 
the following ordinary differential equation:
\begin{equation}
\label{eqn:mpd}
\frac{dY}{dx}
= 
\left( 
\frac{B_{-1}}{(x-a)^2}
+
\frac{B_0}{x-a}
+
\sum_{\nu=1}^3
\frac{A_{\nu}}{x-e_{\nu}} 
\right)
Y(x),
\end{equation}
where
\begin{eqnarray*}
B_{-1} 
&=&
\diag
\left(
\frac{\wp^{\prime}(\alpha)t}{2},
-
\frac{\wp^{\prime}(\alpha)t}{2} 
\right)  \\
B_{0}
&=&
\diag(0,0)\\
A_{\nu} 
&=& 
G^{(\nu)} 
\diag 
\left(-\frac14, \frac14 \right) 
(G^{(\nu)})^{-1}, \,\,
(\nu=1,2,3), \\
G^{(\nu)}
&=&
\left(
\begin{array}{cc}
0 & 
\exp 
\Big\{
\frac{t}{2}
\left(
\frac
{
\wp^{\prime \prime}(\alpha)
}
{2\wp^{\prime}(\alpha)}
+
\zeta(2\alpha)
\right)
\Big\} \\
-
\exp 
\Big\{
-
\frac{t}{2}
\left(
\frac
{
\wp^{\prime \prime}(\alpha)
}
{2\wp^{\prime}(\alpha)}
+
\zeta(2\alpha)
\right)
\Big\} & 0
\end{array}
\right)  \\
& &
\times
\frac
{
\sqrt{2m_{\nu}}
}
{
\sqrt{D^{(\nu)}(\tilde{\omega}_{\nu})i}
}  
\left(
\begin{array}{cc}
\varphi(\tilde{\omega}_{\nu}) &
\varphi(\tilde{\omega}_{\nu})
\left(
\frac{d}{du} 
\log 
\varphi(\tilde{\omega}_{\nu}) 
-
\tilde{\eta_{\nu}}
\right)  \\
\psi(\tilde{\omega}_{\nu}) 
&
\psi(\tilde{\omega}_{\nu})
\left(
\frac{d}{du} 
\log 
\psi(\tilde{\omega}_{\nu}) 
-
\tilde{\eta_{\nu}}
\right)
\end{array}
\right)  
\left(
\begin{array}{cc}
\left(
\frac
{
\wp^{\prime \prime}(\tilde{\omega}_{\nu})
}
{2}
\right)^{\frac14} &  \\
 & 
\left(
\frac
{
\wp^{\prime \prime}(\tilde{\omega}_{\nu})
}
{2}
\right)^{-\frac14}
\end{array}
\right).
\end{eqnarray*}
The behavior near $\infty$ is 
\begin{equation*}
Y(x)
=
G^{(\infty)}
\left(
\left(
\begin{array}{cc}
1 & 0 \\
0 & 1 
\end{array}
\right)
+
O(x^{-1})
\right)
\exp T^{(\infty)}(x),
\end{equation*}
where
\begin{eqnarray*}
G^{(\infty)}
&=&
\left(
\begin{array}{cc}
0
&
\exp
\Big\{
\frac{t}{2}
\left(
\frac
{
\wp^{\prime \prime}(\alpha)
}
{
2\wp^{\prime}(\alpha)
}
+
\zeta(2\alpha)
\right)
\Big\} \\
\exp
\Big\{
-
\frac{t}{2}
\left(
\frac
{
\wp^{\prime \prime}(\alpha)
}
{
2\wp^{\prime}(\alpha)
}
+
\zeta(2\alpha)
\right)
\Big\}
& 0
\end{array}
\right) \\
& &
\times
\frac
{1}
{
\left(
\frac{d}{du}
\log \varphi(0)
-
\
\frac{d}{du}
\log \psi(0)
\right)^{\frac12}
} 
\left(
\begin{array}{cc}
-i\varphi(0) 
& 
i \varphi^{\prime}(0)  \\
-i\psi(0) 
&
i \psi^{\prime}(0)
\end{array}
\right).
\end{eqnarray*}
\end{theorem}

\begin{proof}
From Lemma \ref{lem:ya}, we get
\begin{eqnarray*}
Y^{\prime}(x) Y^{-1}(x)
&=&
\frac{1}{(x-a)^2}
\left(
\begin{array}{cc}
\frac{\wp^{\prime}(\alpha)t}{2} &     \\
                 & -\frac{\wp^{\prime}(\alpha)t}{2} 
\end{array}
\right)
+
\textrm{regular part}  \\
&:=&
\frac{B_{-1}}{(x-a)^2}
+
\frac{B_0}{x-a}
+
\textrm{regular part}.
\end{eqnarray*}
Therefore 
$B_{-1}, B_0$
are determined 
as in Theorem \ref{thm:mpd}
\par
We represent the monodromy matrices 
$M_{\nu} \, (\nu=1, 2, 3)$ 
in the following way:
\begin{eqnarray*}
M_{\nu}
&=&
(C^{(\nu)})^{-1}
\left(
\begin{array}{cc}
\exp \{ -\frac{2\pi i}{4} \}  &     \\
                & \exp \{ \frac{2 \pi i}{4} \}
\end{array}
\right)
C^{(\nu)},
\end{eqnarray*}
where
$C^{(\nu)}$ is the connection matrix.
We set $G^{(\nu)}\,\,(\nu=1,2,3)$ 
in the following way;
\begin{equation}
\label{eqn:nu}
Y(x)
=
G^{(\nu)}
\left(
\left(
\begin{array}{cc}
1 & 0 \\
0 & 1 
\end{array}
\right)
+
O(x-e_{\nu})
\right)
(x-e_{\nu})^
{
\left(
\begin{array}{cc}
-\frac14 &         \\
         & \frac14   
\end{array}
\right)
}
C^{(\nu)}.
\end{equation}
In order to calculate $G^{(\nu)}$, 
we compute the constant term of
\begin{eqnarray*}
& &
Y(x)
(C^{(\nu)})^{-1}
(x-e_{\nu})^
{
\left(
\begin{array}{cc}
\frac14  &          \\
         & -\frac14 
\end{array}
\right)
} 
=
Y(x)
\frac{1}
{
\sqrt{2 i m_{\nu}}
}
\left(
\begin{array}{cc}
m_{\nu} & m_{\nu}  \\
-i         & i  
\end{array}
\right)
(x-e_{\nu})^
{
\left(
\begin{array}{cc}
\frac14  &           \\
         & -\frac14    
\end{array}
\right)
}  \\
&=&
\frac
{
\sqrt{\det \Phi(a)}
}
{
\sqrt{\det \Phi(u)}
}
(G^{(a)})^{-1}  
\frac{1}
{
\sqrt{2 i m_{\nu}}
}
\left(
\begin{array}{cc}
m_{\nu} \varphi(u) -i \varphi(-u) 
& 
m_{\nu} \varphi(u) +i \varphi(-u) \\
m_{\nu} \psi(u)    -i \psi(-u)    
& 
m_{\nu} \psi(u)    +i \psi(-u)
\end{array}
\right) 
(x-e_{\nu})^
{
\left(
\begin{array}{cc}
\frac14 &        \\
        & -\frac14 
\end{array}
\right)
}.
\end{eqnarray*}
\par
By setting 
$
u_{\varphi}=\alpha, u_{\psi}=-\alpha$ in $G^{(a)}
$, 
\begin{equation}
\sqrt{\det \Phi(a)}
\left( G^{(a)} \right)^{-1}
=
\left(
\begin{array}{cc}
0 & 
\exp 
\Big\{
\frac{t}{2}
\left(
\frac
{
\wp^{\prime \prime}(\alpha)
}
{2\wp^{\prime}(\alpha)}
+
\zeta(2\alpha)
\right)
\Big\} \\
-
\exp 
\Big\{
-
\frac{t}{2}
\left(
\frac
{
\wp^{\prime \prime}(\alpha)
}
{2\wp^{\prime}(\alpha)}
+
\zeta(2\alpha)
\right)
\Big\} & 0
\end{array}
\right).
\end{equation}
\par
By direct computation, 
we get 
\begin{eqnarray*}
\det \Phi (u)
&=&
(u-\tilde{\omega}_{\nu})
\left(
\varphi^{\prime}(\tilde{\omega}_{\nu})
\psi(-\tilde{\omega}_{\nu})
-
\varphi(\tilde{\omega}_{\nu})
\psi^{\prime}(-\tilde{\omega}_{\nu})    
+
\varphi^{\prime}(-\tilde{\omega}_{\nu})
\psi(\tilde{\omega}_{\nu})
-
\varphi(-\tilde{\omega}_{\nu})
\psi^{\prime}(\tilde{\omega}_{\nu})
\right)
+
\cdots \\
&=&
(u-\tilde{\omega}_{\nu})
\Big\{
\varphi(\tilde{\omega}_{\nu})
\psi(-\tilde{\omega}_{\nu})
\left(
\frac
{
\varphi^{\prime}(\tilde{\omega}_{\nu})
}
{
\varphi(\tilde{\omega}_{\nu})
} 
-
\frac
{
\psi^{\prime}(-\tilde{\omega}_{\nu})
}
{
\psi(-\tilde{\omega}_{\nu})
}
\right)   \\
& &
\hspace{30mm}
+
\varphi(-\tilde{\omega}_{\nu})
\psi(\tilde{\omega}_{\nu})
\left(
\frac
{
\varphi^{\prime}(-\tilde{\omega}_{\nu})
}
{
\varphi(-\tilde{\omega}_{\nu})
}
-
\frac
{
\psi^{\prime}(\tilde{\omega}_{\nu})
}
{
\psi(\tilde{\omega}_{\nu})
}
\right)
\Big\}
+
\cdots.
\end{eqnarray*}
From 
(\ref{eqn:per1}), (\ref{eqn:per2}), 
\begin{eqnarray*}
\det \Phi(u)
&=&
(u-\tilde{\omega}_{\nu})
\frac
{
m_{\nu}
}
{m_{\infty}}
\varphi(\tilde{\omega}_{\nu})
\psi(\tilde{\omega}_{\nu})
\times
2
\left(
\frac{d}{du}
\log
\varphi(\tilde{\omega}_{\nu})
-
\frac{d}{du}
\log
\psi(\tilde{\omega}_{\nu})
\right) 
+ 
O((u-\tilde{\omega}_{\nu})^2) \\
&=&
(u-\tilde{\omega}_{\nu})
D^{(\nu)}(\tilde{\omega}_{\nu})  
+ 
O((u-\tilde{e_{\nu}})^2)   \\
&=&
(x-e_{\nu})^{\frac12}
\left(
\frac
{
\wp^{\prime \prime}
(
\tilde{\omega}_{\nu}
)
}
{2}
\right)^{-\frac12}
D^{(\nu)}(\tilde{\omega}_{\nu})
+ 
O((\sqrt{x-e_{\nu}})^2)
\end{eqnarray*}
From 
(\ref{eqn:per1}), (\ref{eqn:per2}), 
we note that
\begin{eqnarray*}
m_{\nu} 
\varphi(\tilde{\omega}_{\nu}) 
+ 
i 
\varphi(-\tilde{\omega}_{\nu}) 
&=& 
0, \\
m_{\nu} 
\psi(\tilde{\omega}_{\nu}) 
+ 
i 
\psi(-\tilde{\omega}_{\nu})       
&=& 
0.
\end{eqnarray*}
Therefore, 
\begin{eqnarray*}
G^{(\nu)}
&=&
\left.
Y(x)
(C^{(\nu)})^{-1}
(x-e_{\nu})^
{
\left(
\begin{array}{cc}
\frac14  &     \\
         & -\frac14 
\end{array}
\right)
}
\right|_{x=e_{\nu}} \\
&=&
(\det \Phi(a))^{\frac12}
(G^{(a)})^{-1}
(
D^{(\nu)}(\tilde{\omega}_{\nu})
)^{-\frac12}
(
\frac
{
\wp^{\prime \prime}(\tilde{\omega}_{nu})
}
{2}
)^{\frac14}  \\
& &
\hspace{10mm}
\times
\frac{1}{\sqrt{2i m_{\nu}}}
\left(
\begin{array}{cc}
m_{\nu} \varphi(u)-i \varphi(-u) &
\frac
{
m_{\nu} \varphi(u) + i \varphi (-u)
}
{
u-\tilde{\omega}_{\nu}
}  \\
m_{\nu} \psi(u)-i \psi(-u) &
\frac
{
m_{\nu} \psi(u) + i \psi (-u)
}
{
u-\tilde{\omega}_{\nu}
} 
\end{array}
\right)
\left.
\left(
\begin{array}{cc}
1   &    \\
    & 
\frac{
u-\tilde{\omega}_{\nu}
}
{
(x-e_{\nu})^{\frac12}
}
\end{array}
\right)
\right|_{x=e_{\nu}} \\
&=&
(\det(a))^{\frac12}
(G^{(a)})^{-1}
\frac
{
\sqrt{2m_{\nu}}
}
{
\sqrt{D^{(\nu)}(\tilde{\omega}_{\nu})i}
}
\left(
\begin{array}{cc}
\varphi(\tilde{\omega}_{\nu})  & 
\varphi(\tilde{\omega}_{\nu})
\left(
\frac{d}{du}
\log \varphi(\tilde{\omega}_{\nu})
-
\tilde{\eta_{\nu}}
\right)  \\
\psi(\tilde{\omega}_{\nu})  & 
\psi(\tilde{\omega}_{\nu})
\left(
\frac{d}{du}
\log \psi(\tilde{\omega}_{\nu})
-
\tilde{\eta_{\nu}}
\right)
\end{array}
\right)  \\
& &
\hspace{60mm}
\times
\left(
\begin{array}{cc}
\left(
\frac{\wp^{\prime \prime}(\tilde{\omega}_{\nu})}{2}
\right)^{\frac14} &    \\
 & \left(
\frac{\wp^{\prime \prime}(\tilde{\omega}_{\nu})}{2}
\right)^{-\frac14}
\end{array}    
\right).
\end{eqnarray*}
From (\ref{eqn:nu}), we get
\begin{eqnarray*}
Y^{\prime}(x)Y^{-1}(x) 
&=&
\frac{1}{(x-e_{\nu})}
G^{(\nu)}
\left(
\begin{array}{cc}
-\frac14 &          \\
        & \frac14
\end{array}
\right)
(G^{(\nu)})^{-1}
+ 
\textrm{regular part}   \\
&:=&
\frac{A_{\nu}}{(x-e_{\nu})}
+
\textrm{regular part}.
\end{eqnarray*}
\par
We set $G^{(\infty)}$ in the following way:
\begin{equation*}
Y(x)
=
G^{(\infty)}
\left(
1
+
O(x^{-1})
\right)
\exp T^{(\infty)}(x).
\end{equation*}
In the same way, we can compute $G^{(\infty)}$.
\end{proof}

We can get the following deformation equation 
from Theorem \ref{thm:mpd}.

\begin{corollary}
The deformation equation 
of the monodromy preserving deformation (\ref{eqn:mpd}) 
is as follows. For $\nu=1,2,3,$
\begin{eqnarray}
d A_{\nu}
&=&
\sum_{\mu \neq \nu}
\frac
{
[A_{\mu},A_{\nu}]
}
{e_{\nu}-e_{\mu}}
de_{\nu} 
+
\sum_{\mu}
\frac
{
[A_{\mu},A_{\nu}]
}
{a-e_{\mu}}
de_{\mu} \notag \\
&-&
\frac
{
[A_{\nu},B_{-1}]
}
{
(a-e_{\nu})^2
}
de_{\nu}
+
\frac
{
[d T^{(a)}_{-1}, A_{\nu}]
}
{a-e_{\nu}} 
+
[
[d T^{(a)}_{-1}, Y^{(a)}_{1}], 
A_{\nu}
],
\end{eqnarray}
where $d$ is the exterior differentiation 
with respect to the deformation parameters, 
$t, e_1, e_2, e_3$ 
and 
$T^{(a)}_{-1}, Y^{(a)}_1$ is 
defined in  Lemma \ref{lem:ya}.
\end{corollary}

\section{The  tau-function for the Schlesinger System}
In this section, 
we will calculate 
the $\tau$-function for 
(\ref{eqn:mpd}). 
This section 
consists of 
three subsections. 
Subsection 6.1 is devoted to 
the Hamiltonian $H_t$. 
Subsection 6.2, 6.3 
is devoted to the Hamiltonian 
$H_{\nu} \, (\nu=1,2,3).$ 
In subsection 6.2, 
we show some facts about elliptic functions. 
In subsection 6.3, 
we calculate $H_{\nu}$ and the $\tau$-function.

\subsection{The Hamiltonian at the Irregular Singular Point}
In this subsection, 
we compute $\omega_a, H_t$ 
in the following way:
\begin{proposition}
\label{prop:omegaa}
\begin{eqnarray*}
\omega_a &=& 
\frac
{
\sigma [p,q]^{\prime}(t)
}
{\sigma [p,q] (t)} dt
+
\frac
{
\sigma [p,q]^{\prime}(t) 
}
{\sigma [p,q](t)} 
t
\left(
-
\frac{d e_1}{2(a-e_1)} 
-
\frac{d e_2}{2(a-e_2)} 
-
\frac{d e_3}{2(a-e_3)}
\right) \\
&+&
f(e_1, e_2, e_3)
\left(
\frac{t dt}{2} 
-
\frac{t^2 d e_1}{4(a-e_1)}
-
\frac{t^2 d e_2}{4(a-e_2)}
-
\frac{t^2 d e_3}{4(a-e_3)}
\right),
\end{eqnarray*}
where
\begin{equation}
f(e_1, e_2, e_3)
=
-
a 
+
\frac13 
\sum_{\nu=1}^3 e_{\nu}
+
\frac12 
\prod_{\nu=1}^3 (a-e_{\nu}) 
\sum_{\mu=1}^3 \frac{1}{(a-e_{\mu})^2}.
\end{equation}
From $\omega_a$, we get
\begin{eqnarray*}
H_t
&=&
\frac
{
\sigma [p,q]^{\prime}(t)
}
{\sigma [p,q](t)} 
+ 
\frac{t}{2} f(e_1, e_2, e_3) \\
&=&
\frac{\partial}{\partial t}
\left\{
\frac{\eta_1 t^2}{2 \omega_1}
\right\}
+
\frac{\partial}{\partial t}
\Big\{
\log
\theta[p,q]
\left(
\frac{t}{\omega_1}
\right)
\Big\} 
+
\frac{\partial}{\partial t}
\Big\{
\frac{t^2}{4}
f(e_1,e_2,e_3)
\Big\}.
\end{eqnarray*}
\end{proposition}

\begin{proof}
From the definition of $ \omega_a $,
\begin{equation}
\omega_a
= 
- 
\Res_{x=a} \,\, 
\trace \hat{Y}^{(a)}(x)^{-1} 
\frac
{
\partial \hat{Y}^{(a)}
}
{\partial x}(x)
\, d T^{(a)}(x).
\end{equation}
We set 
\begin{equation}
\Pi (P) := -\frac{t}{2}
\{ 
\zeta(u-\alpha)
+ 
\zeta(u+\alpha) 
\}.
\end{equation}
From Lemma \ref{lem:inv}, we get 
\begin{eqnarray*}
\Pi (P)
&=& 
-\frac{t}{2}
\{
\frac{1}{u-\alpha}
+\cdots
+
\zeta(u+\alpha)
\} \\
&=&
- 
\frac
{
\wp^{\prime}(\alpha)t
}
{2} 
\frac{1}{(x-a)}
- 
\frac{t}{2} 
\left( 
\frac
{
\wp^{\prime \prime}(\alpha)
}
{2 \wp^{\prime}(\alpha)}
+ 
\zeta(2\alpha) 
\right) \\
& & 
- 
\frac{t}{2} 
\left( 
-
\frac
{
\wp(2\alpha)
}
{\wp^{\prime}(\alpha)}
-
\frac
{
\wp^{\prime \prime}(\alpha)^2
}
{
4 \wp^{\prime}(\alpha)^3
}
+
\frac{
\wp^{\prime \prime \prime}(\alpha)
}
{
6 \wp^{\prime}(\alpha)^2
}
\right) 
(x-a) 
+
\cdots.
\end{eqnarray*}
We denote the regular part of $\Pi(P)$ 
around $x=a$ by $\hat{\Pi}(P)$:
\begin{equation}
\hat{\Pi}(P) 
:= 
\Pi(P) 
- 
\left(
- 
\frac{\wp^{\prime}(\alpha)t}{2} 
\frac{1}{(x-a)} 
\right).
\end{equation} 
We set 
\begin{eqnarray*}
\hat{\varphi}(P)
&=&
\sigma[p,q](u+u_{\varphi}+t)
\sigma(u-u_{\varphi}), \\
\hat{\psi}(P)
&=&
\sigma[p,q](u+u_{\psi}+t),
\sigma(u-u_{\psi}).
\end{eqnarray*}
Then we get
\begin{eqnarray*}
Y(P) 
&=&
\frac
{
\sqrt{\det \Phi(a)}
}
{
\sqrt{\det \Phi(P)}
} 
(G^{(a)})^{-1} 
\Phi(P) \\
&=&
\frac
{
\sqrt{\det \Phi(a)}
}
{
\sqrt{\det \Phi(P)}
} 
(G^{(a)})^{-1}
\left(
\begin{array}{cc}
\hat{\varphi}(P) \exp (\hat{\Pi} (P)) 
& 
\hat{\varphi}(P^{*}) \exp (\hat{\Pi} (P^{*})) \\
\hat{\psi}(P) \exp (\hat{\Pi} (P)) 
& 
\hat{\psi}(P^{*}) \exp (\hat{\Pi} (P^{*}))
\end{array}
\right) \\
& &         
\diag 
\left(
\exp
\Big\{
- 
\frac
{
\wp^{\prime}(\alpha)t
}
{2} 
\frac{1}{x-a}
\Big\}, 
\exp
\Big\{
\frac
{
\wp^{\prime}(\alpha)t
}
{2} 
\frac{1}{x-a}
\Big\}
\right) \\
&:=&
\hat{Y}^{(a)}(x) 
\exp T^{(a)} (x).
\end{eqnarray*}
According to the definition of $\omega_a,$ 
\begin{equation}
\omega_a 
= 
- 
\Res_{x=a} 
\trace \hat{Y}^{(a)}(x)^{-1} 
\frac{\partial}{\partial x} 
\hat{Y}^{(a)}(x) 
d T^{(a)} (x).
\end{equation}
We set 
\begin{equation}
A(x) 
= 
(G^{(a)})^{-1}
\left(
\begin{array}{cc}
\hat{\varphi}(u) \exp \hat{\Pi} (u) 
& 
\hat{\varphi}(-u) \exp \hat{\Pi}(-u) \\
\hat{\phi} (u) \exp \hat{\Pi} (u)   
& \hat{\phi}(-u) \exp \hat{\Pi}(-u)
\end{array}
\right).
\end{equation}
By direct calculation, we get
\begin{equation}
- 
\trace 
\hat{Y}^{(a)}(x)^{-1} 
\frac{\partial}{\partial x} 
\hat{Y}^{(a)} (x) 
d T^{(a)} (x)  
=
\frac
{
d 
\left(
\wp^{\prime}(\alpha) t 
\right)
}
{2}  
\frac{1}{(x-a)} 
\trace A^{-1}(x) A^{\prime}(x) 
\left(
\begin{array}{cc}
1 &  \\
  & -1 
\end{array}
\right).
\end{equation}
We obtain
\begin{eqnarray*}
& &  
\trace A^{-1}(x) A^{\prime}(x) 
\left(
\begin{array}{cc}
1 &  \\
  & -1 
\end{array}
\right)  \\
&=& 
\frac{1}{\det \Phi (u)}
\Big[
\det
\left(
\begin{array}{cc}
\{ 
\hat{\varphi} (u) \exp \hat{\Pi}(u) 
\}^{\prime} 
& \hat{\varphi}(-u) \exp \hat{\Pi}(-u) \\
\{ 
\hat{\psi} (u) \exp \hat{\Pi}(u) 
\}^{\prime}    
& 
\hat{\psi}(-u) \exp \hat{\Pi}(-u) 
\end{array}
\right) \\      
& &  
\hspace{18mm}           
-
\det
\left(
\begin{array}{cc}
\hat{\varphi} (u) \exp \hat{\Pi}(u) 
& 
\{ 
\hat{\varphi}(-u) \exp \hat{\Pi}(-u) 
\}^{\prime} \\
\hat{\psi} (u) \exp \hat{\Pi} (u) 
& 
\{ \hat{\psi}(-u) \exp \hat{\Pi} (-u) \}^{\prime}
\end{array}
\right)
\Big],
\end{eqnarray*}
where $\prime$ means the differentiation 
with respect to the variable $x$.
\par
We normalized the matrix function $Y(x)$ around $x=a$ 
like Lemma \ref{lem:ya} and proved that 
the monodromy group of $Y(x)$ 
is independent of $P_{\varphi}, P_{\psi}$ 
in Theorem \ref{thm:mono} and its corollary. 
Therefore, 
we can choose the parameters, 
$P_{\varphi}, P_{\psi}$ 
at our disposal to simplify the calculation. 
Firstly, we multiply 
both the numerators and the denominators 
by $\frac{1}{x_{\psi}-x_{\varphi}}$. 
We take the limit $P_{\psi} \rightarrow P_{\varphi}$ 
and get
\begin{equation}
\hat{\psi}(P) 
= 
\frac
{
\partial \hat{\varphi}(P)
}
{\partial x_{\varphi}}.
\end{equation}
Secondly, we multiply 
both 
the numerators 
and 
the denominators 
by $\frac{1}{x_{\varphi}-x}$. 
We take the limit 
$
P_{\varphi} \rightarrow P
$ 
and obtain
\begin{equation}
\trace A^{-1}(x) A^{\prime}(x)
\left(
\begin{array}{cc}
1 &  \\
  & -1 
\end{array}
\right)
=
2 
\frac{1}{\sigma[p,q](t)} 
\lim_{P_{\varphi} \rightarrow P} 
\frac{\partial}{\partial x_{\varphi}}
\sigma[p,q](-u + u_{\varphi} + t) 
+ 
2 
\frac{\partial}{\partial x}
\{
\hat{\Pi} (P)
\}.
\end{equation}
According to the definition of $\omega_a$,
\begin{eqnarray*}
\omega_a 
&=& 
\frac{d (\wp^{\prime}(\alpha) t)}{2}
\trace A^{-1}(x) A^{\prime}(x) 
\left.
\left(
\begin{array}{cc}
1 &  \\ 
  &  -1
\end{array}  
\right) 
\right|_{x=a}  \\
&=&
\frac
{
\sigma[p,q]^{\prime}(t)
}
{\sigma[p,q](t)} dt 
+
\frac
{
\sigma[p,q]^{\prime}(t)
}
{\sigma[p,q](t)} t 
\frac
{
d \wp^{\prime}(\alpha)
}
{\wp^{\prime}(\alpha)}  \\
& &  
\hspace{15mm}     
+ 
\left(
\frac{t}{2} dt 
+
\frac{t^2}{2}
\frac{d \wp^{\prime}(\alpha)}
{\wp(\alpha)}
\right)
\left(
\wp(2\alpha) 
+ 
\frac
{
\wp^{\prime \prime}(\alpha)^2
}
{4\wp^{\prime}(\alpha)^2} -
\frac
{
\wp^{\prime \prime \prime}(\alpha)
}
{6 \wp^{\prime}(\alpha)}
\right) 
\end{eqnarray*}
where $d$ is the exterior differentiation 
with respect to the deformation parameters, 
$e_1, e_2, e_3, t$. 
\par
In order to compute $\omega_a$, 
we need some preparations about 
the $\wp$-function. 
By direct calculation, we get
\begin{equation}
\label{eqn:ex}
\frac
{d \wp^{\prime}(\alpha)
}
{\wp^{\prime}(\alpha)}
=
-\frac12 
\frac{1}{a-e_1}
d e_1
-
\frac12 
\frac{1}{a-e_2}
d e_2
-
\frac12 
\frac{1}{a-e_3}
d e_3.
\end{equation}
We will use the addition theorem
\begin{equation}
\label{eqn:add}
\wp(2\alpha) =- 2\wp(\alpha) + \frac14 
\left(
\frac
{
\wp^{\prime \prime}(\alpha)
}
{\wp^{\prime}(\alpha)}
\right)^2,
\end{equation}
and the following equation:
\begin{equation}
\label{eqn:tri}
\wp^{\prime \prime \prime}(\alpha)
= 
12 
\wp^{\prime}(\alpha) 
\wp(\alpha).
\end{equation}
From the relationship 
between 
the $\wp$-function and the Abel-map, 
we obtain
\begin{eqnarray*}
\wp(\alpha) 
&=&
a
- 
\frac13 \sum_{i=1}^3 e_i  \\
\wp^{\prime}(\alpha)^2 
&=& 
4(a-e_1)(a-e_2)(a-e_3) \\
\wp^{\prime \prime} (\alpha) 
&=& \sum_{i<j} 2(a-e_i)(a-e_j).
\end{eqnarray*}
\par
By using the above equations, we get
\begin{eqnarray*}
\omega_a &=& 
\frac
{
\sigma [p,q]^{\prime}(t)
}
{\sigma [p,q] (t)} 
dt
+
\frac
{
\sigma [p,q]^{\prime}(t)
}
{\sigma [p,q](t)}
t
\left(
-
\frac{d e_1}{2(a-e_1)} 
-
\frac{d e_2}{2(a-e_2)} 
-
\frac{d e_3}{2(a-e_3)}
\right) \\
&+&
f(e_1, e_2, e_3)
\left(
\frac{t dt}{2} 
-
\frac{t^2 d e_1}{4(a-e_1)}
-
\frac{t^2 d e_2}{4(a-e_2)}
-
\frac{t^2 d e_3}{4(a-e_3)}
\right).
\end{eqnarray*}
\end{proof}

\subsection{Three Lemmas about Elliptic Functions}
We show 
three lemmas about elliptic functions.

\begin{lemma}
\label{lem:sig}
\begin{equation}
\label{eqn:u3}
\omega_1 \eta_1 
= 
- 
\frac13 
\frac
{
\theta_{11}^{\prime \prime \prime}
}
{
\theta_{11}^{\prime}
}
\end{equation}
\begin{equation}
\label{eqn:u5}
- 
\frac13 
\left(
\sum_{i=1}^3 
e_i 
\right)^2 
+ 
(e_1 e_2 +e_2 e_3 + e_3 e_1)
= 
\frac12 
\left( 
\frac
{
\theta_{11}^{(5)}
}
{\theta_{11}^{(1)}} 
\right)
- 
\frac56 
\left( 
\frac
{
\theta_{11}^{(3)}
}
{\theta_{11}^{(1)}} 
\right)^2 
\frac{1}{\omega_1^4}.
\end{equation}
\end{lemma}

\begin{proof}
We can show 
the first equation (\ref{eqn:u3}) 
by comparing the coefficients on $u^3$ 
of (\ref{eqn:sig}).
\par
The second equation (\ref{eqn:u5}) 
can 
be proved 
by comparing the coefficients on $u^5$ 
of (\ref{eqn:sig}) 
and 
by using the following equation:
\begin{equation}
\tilde{e_1}\tilde{e_2}
+
\tilde{e_2}\tilde{e_3}
+
\tilde{e_3}\tilde{e_1}
= 
-
\frac13 
\left(
\sum_{i=1}^3 e_i 
\right)^2 
+ 
\left(
e_1 e_2 +e_2 e_3 +e_3 e_1 
\right).
\end{equation}
\end{proof}

\begin{lemma}
For $\nu=1, 2, 3,$ the dependence of 
the period $\Omega$ on the branch points is 
described by the equation
\begin{equation}
\label{eqn:Omega}
\frac{\partial \Omega}{\partial e_{\nu}} 
=
\frac
{
\pi i
}
{
\omega_1^2 
\prod_{\mu \neq \nu} 
(e_{\nu} - e_{\mu})
}.
\end{equation}
\end{lemma}

\begin{proof}
We will show the case $\nu=1$. 
The other cases can be proved 
in the same way.
\par
According to the definition of the Abel map, 
we have 
\begin{equation}
du(x) = 
\frac{dx}{\sqrt{4(e-e_1)(x-e_2)(x-e_3)}}.
\end{equation}
The dependence of $du$ on $e_1$ is 
\begin{equation}
\frac{\partial}{\partial e_1} \{du(x) \}=
\frac{1}{2(x-e_1)} du(x).
\end{equation}
Now, we consider the integral
\begin{equation}
\oint_{\partial \hat{E}} u(x) 
\frac{\partial}{\partial e_1} du
=
\oint_{\partial \hat{E}} 
\frac{1}{2(x-e_1)} u(x) du(x).
\end{equation}
The left side is
\begin{equation}
-
\omega_2 
\frac{\partial \omega_1}{\partial e_1} 
+ 
w_1
\frac
{
\partial \omega_2
}
{\partial e_1}.
\end{equation} 
\par
We set 
the local coordinate of the branch point $e_1$ 
in the following way:
\begin{equation}
t_1 = \sqrt{x-e_1}.
\end{equation}
The right side is 
by the residue theorem 
\begin{equation}
\pi i 
\frac{du}{dt_1}(e_1) 
\frac{du}{dt_1}(e_1)
= 
\frac
{
\pi i \omega_1^2
}
{
\omega_1^2 (e_1-e_2)(e_1-e_3)
}.
\end{equation}
\par 
By comparing the both sides, we get
\begin{equation}
- 
\frac{\omega_2}{\omega_1^2} 
\frac{\partial \omega_1}{\partial e_1}
+ 
\frac{\partial \omega_2}{\partial e_1} 
\frac{1}{\omega_1}
=
\frac{\pi i }{\omega_1^2 (e_1-e_2)(e_1-e_3)}.
\end{equation}
Then we obtain
\begin{equation}
\frac{\partial \Omega}{\partial e_1}
=
\frac{\pi i}{\omega_1^2 (e_1-e_2)(e_1-e_3)}.
\end{equation}
\end{proof}

\begin{lemma}
\label{lem:frac}
For $\nu=1,2,3,$ 
\begin{equation}
\frac{\partial }{\partial e_\nu}
\left(
\frac{\eta_1 t^2}{2 \omega_1}
\right)
=
t^2 \left(
\frac{\partial}{\partial e_{\nu}}
\log 
\omega_1
\right)^2
\prod_{\mu \neq \nu} (e_{\nu}-e_{\mu})
- \frac{t^2}{12}.
\end{equation}
\end{lemma}

\begin{proof}
We will show the case $\nu=1$. 
The other cases can be proved 
in the same way.
From (\ref{eqn:u3}), (\ref{eqn:u5}) 
and the heat equation
\begin{equation}
\label{eqn:heat}
\frac
{
\partial \theta [p,q]
}
{\partial z^2}
(z; \Omega)
=
4 \pi i
\frac
{
\partial \theta [p,q]
}
{\partial \Omega} 
(z; \Omega),
\end{equation}
we get
\begin{eqnarray*}
\frac{\partial }{\partial e_1}
\left(
\frac{\eta_1 t^2}{2 \omega_1}
\right)
&=&
\frac{t^2}{3}
\frac
{
\frac{\partial \omega_1}{\partial e_1}
}
{\omega_1^3}
\frac
{
4 \pi i 
\frac
{
\partial \theta_{11}^{\prime}
}
{\partial \Omega}
}
{\theta_{11}^{\prime}}
-
\frac{t^2}{6 \omega_1^2}
\frac{\partial \Omega}{\partial e_1}
\frac{1}{4 \pi i}
\Big\{
2 \omega_1^4 \left(- \frac{g_2}{4} \right)
+
\frac23 
\left( 
\frac
{
\theta_{11}^{\prime \prime \prime}
}
{\theta_{11}^{\prime}} 
\right)^2
\Big\}. \\
\end{eqnarray*}
By differentiating 
the logarithm of 
both sides of the formula:
\begin{equation*}
\sqrt[4]
{
16
(e_1-e_2)^2 (e_1-e_3)^2 (e_2- e_3)^2
}
=
\frac{2 \pi}{\omega_1^3} 
(\theta_{11}^{\prime})^2,
\end{equation*}
we get
\begin{equation}
\label{eqn:log}
3 
\frac
{
\frac{\partial \omega_1}{\partial e_1}
}
{\omega_1}
+
\frac12 
\frac{1}{(e_1-e_2)} 
+ 
\frac12 
\frac{1}{(e_1-e_3)}
=
2
\frac
{
\frac
{
\partial \theta_{11}^{\prime}
}
{\partial e_1}
}
{\theta_{11}^{\prime}}.
\end{equation}
\par
By using 
(\ref{eqn:Omega}) and (\ref{eqn:log}), 
we obtain
\begin{eqnarray*}
\frac{\partial }{\partial e_1}
\left(
\frac{\eta_1 t^2}{2 \omega_1}
\right)
&=&
t^2 
\left(
\frac
{
\frac{\partial \omega_1}{\partial e_1}
}
{\omega_1}
\right)^2
(e_1-e_2)(e_1-e_3)
-\frac{t^2}{12}.
\end{eqnarray*}
\end{proof}

\subsection{The tau-function}
In order to compute $\omega_{e_{\nu}} (\nu=1,2,3)$, 
we will show the following lemma.

\begin{lemma}
\label{lem:omeganu}
For $\nu=1,2,3,$ 
\begin{eqnarray*}
\Res_{x=e_{\nu}}
\,
\frac12
\left(
\frac{d Y}{dx}
Y^{-1}(x)
\right)^2
&=&
-
\frac18
\left(
\sum_{\mu \neq \nu}
\frac{1}{e_{\nu}-e_{\mu}}
\right)
-
\frac12
\frac{\partial}{\partial e_{\nu}}
\log \omega_1  \\
& &
+
t^2
\left(
\frac{\partial}{\partial e_{\nu}}
\log \omega_1
\right)^2
\prod_{\mu \neq \nu}
(e_{\nu}-e_{\mu})
+
\frac{\partial}{\partial e_{\nu}}
\big\{
\log
\theta[p,q]
\left(
\frac{t}{\omega_1};
\Omega
\right)
\big\}  \\
& &
+
\frac{t}{2(a-e_{\nu})}
\frac
{
\sigma[p,q]^{\prime}(\frac{t}{\omega_1};\Omega)
}
{
\sigma[p,q](\frac{t}{\omega_1};\Omega)
}
+
\frac{t^2}{4}
\frac
{
\prod_{\mu \neq \nu}
(e_{\nu}-e_{\mu})
}
{(a-e_{\nu})^2} \\
& &
+
\frac{t^2}{6(a-e_{\nu})}
\left(
\sum_{\mu \neq \nu}
(e_{\nu}-e_{\mu})
\right).
\end{eqnarray*}
\end{lemma}

\begin{proof}
We will prove the case $\nu=1$. 
The other cases can be shown in the same way.
By direct calculation, we obtain 
\begin{equation*}
\frac12 
\trace \left( Y^{\prime}(x)Y^{-1}(x) \right)^2
=
-
\frac{\det \left( \Phi_x \right)}{\det \Phi}
+
\frac14
\left(
\frac
{
(\det \Phi)_x
}
{\det \Phi}
\right)^2.
\end{equation*}
We will calculate $\omega_1$ just like $\omega_a$ 
in Proposition \ref{prop:omegaa}.
We multiply 
both the numerators and the denominators of
$$
\frac{\det \left( \Phi_x \right)}{\det \Phi},
\frac
{
(\det \Phi)_x
}
{\det \Phi}
$$
by $\frac{1}{x_{\varphi}-x_{\psi}}$. 
We take the limit $P_{\psi} \rightarrow P_{\varphi}$ 
and get
\begin{equation}
\psi(P) 
= 
\frac{\partial \varphi(P)}{\partial x_{\varphi}}.
\end{equation}
We take the limit $P_{\varphi} \rightarrow P$ 
and obtain 
\begin{eqnarray*}
\frac
{
(\det \Phi)_x
}
{\det \Phi}
&=&
2 
\frac{1}{\sigma(-2u)}
\lim_{P_{\varphi} \rightarrow P}
\frac{\partial}{\partial x}
\sigma (-u-u_{\varphi})  \\
\frac{\det \left( \Phi_x \right)}{\det \Phi} 
&=& 
\frac{1}{\sigma[p,q](t)}
\lim_{P_{\varphi} \rightarrow P}
\frac{\partial^2}{\partial x \partial x_{\varphi}}
\sigma[p,q](-u+u_{\varphi}+t)  \\
&+& 
\frac{2}{\sigma[p,q](t)} 
\lim_{P_{\varphi} \rightarrow P}
\frac{\partial }{\partial x}
\sigma[p,q] (-u+u_{\varphi}+t)
\frac{\partial }{\partial x} \Pi (P) \\
&+&
\frac{1}{\sigma(-2u)}
\lim_{P_{\varphi} \rightarrow P}
\frac{\partial^2}{\partial x \partial x_{\varphi}}
\sigma (-u-u_{\varphi}) \\
&-&
\left(
\frac{\partial}{\partial x} \Pi (P)
\right)^2.
\end{eqnarray*}
Therefore, we get
\begin{eqnarray*}
\frac12 
\trace \left( Y^{\prime}(x)Y^{-1}(x) \right)^2
&=&
-
\lim_{P_{\varphi} \rightarrow P}
\frac{\partial^2}{\partial x \partial x_{\varphi}}
\log \sigma (-u -u_{\varphi}) \\
& & -
\frac{1}{\sigma[p,q](t)}
\lim_{P_{\varphi} \rightarrow P}
\frac{\partial^2}{\partial x \partial x_{\varphi}}
\sigma[p,q] (-u + u_{\varphi} +t)  \\
& & -
\frac{2}{\sigma[p,q](t)}
\lim_{P_{\varphi} \rightarrow P}
\frac{\partial}{\partial x} \sigma[p,q](-u+u_{\varphi}+t)
\frac{\partial}{\partial x}
\Pi (P)  \\
& & +
\left(
\frac{\partial}{\partial x}
\Pi (P)
\right)^2.
\end{eqnarray*}
In order to calculate the first term 
of 
$
\frac12 
\trace \left( Y^{\prime}(x)Y^{-1}(x) \right)^2,
$ 
we will use the following formulas:
\begin{eqnarray*}
\wp(u) 
&=& 
- 
\frac{d^2}{d u^2} \log \sigma(u)  \\
\wp(u) 
&=& 
x 
- 
\frac13 
\sum_{\mu=1}^3 e_{\mu}, 
\,\, \tilde{e_{\nu}}=e_{\nu}
-\frac13 
\sum_{\mu=1}^3 e_{\mu}
\,\, (\nu=1,2,3)  \\
\wp^{\prime}(u)^2 
&=& 
4 
(\wp(u)-\tilde{e_1})
(\wp(u)-\tilde{e_2})
(\wp(u)-\tilde{e_3})
\end{eqnarray*}
Then 
\begin{eqnarray}
\label{eqn:A}
& &
\lim_{P_{\varphi} \rightarrow P}
\frac{\partial^2}{\partial x \partial x_{\varphi}}
\log \sigma (-u -u_{\varphi}) \notag \\
&=&
\left(
\frac{du}{dx}
\right)^2
\frac
{
\sigma^{\prime \prime}(-2u) 
\sigma(-2u) 
- 
\sigma^{\prime}(-2u)^2
}
{\sigma(-2u)^2} \notag \\
&=&
\frac{1}{24}
\sum_{\nu=1}^3
\frac{1}{(x-e_{\nu})}
\sum_{\mu \neq \nu} 
\frac{1}{(e_{\nu}-e_{\mu})}
-
\frac{1}{16} 
\sum_{\nu=1}^3 
\frac{1}{(x-e_{\nu})^2}. 
\end{eqnarray}
\par
By direct calculation, 
we can calculate the second term of 
$
\frac12 
\trace 
\left( 
Y^{\prime}(x)Y^{-1}(x) 
\right)^2,
$ 
in the following way:
\begin{eqnarray}
\label{eqn:B}
& & \frac{1}{\sigma[p,q](t)}
\lim_{P_{\varphi} \rightarrow P}
\frac
{
\partial^2
}
{
\partial x_{\varphi} 
\partial x
}
\sigma[p,q](-u+u_{\varphi}+t)  \notag \\
& &
\hspace{40mm}
=
-
\frac
{
\sigma[p,q]^{\prime \prime}(t)
}
{\sigma[p,q](t)}
\left(
\sum_{\nu=1}^3
\frac{1}{x-e_{\nu}}
\frac{1}
{
4\prod_{\mu \neq \nu }(e_{\nu}-e_{\mu})
}
\right).
\end{eqnarray}
\par
In order to calculate the third and forth terms of
$
\frac12 
\trace 
\left( 
Y^{\prime}(x)Y^{-1}(x) 
\right)^2,
$ 
we use 
the following relation:
\begin{equation}
\label{eqn:zeta}
\wp(u) 
= 
- 
\frac{d}{du} 
\zeta(u).
\end{equation}
By using (\ref{eqn:zeta}), 
we get
\begin{eqnarray}
\label{eqn:C}
& &
\frac{2}{\sigma[p,q](t)}
\lim_{P_{\varphi} \rightarrow P}
\frac{\partial}{\partial x}
\sigma[p,q](-u+u_{\varphi}+t)
\frac{\partial}{\partial x}
\Pi (P)  \notag \\
& & =
t
\frac
{
\sigma[p,q]^{\prime}(t)
}
{\sigma[p,q](t)}
\Big\{
\frac{x}
{
2 \prod_{\nu=1}^3 (x-e_{\nu})
}
+
\frac{a}{2 \prod_{\nu=1}^3 (x-e_{\nu})}
-
\frac13
\frac
{
\sum_{\nu=1}^3 e_{\nu}
}
{\prod_{\nu=1}^3 (x-e_{\nu})}  \notag \\
& &
\hspace{40mm}
-
\frac{1}{2(x-a)^2}
-
\frac12
\frac{1}{(x-a)^2}
\frac
{
\prod_{\nu=1}^3 (a-e_{\nu})
}
{\prod_{\nu=1}^3(x-e_{\nu})}
\Big\}
\end{eqnarray}
and
\begin{eqnarray}
\label{eqn:D}
\Res_{x=e_{1}}
\left(
\frac{\partial}{\partial x}
\Pi (P)
\right)^2
&=&
\frac{t^2}{36}
\frac{e_1-e_2}{e_1-e_3}
+
\frac{t^2}{36}
\frac{e_1-e_3}{e_1-e_2}
+
\frac{t^2}{4}
\frac{(e_1-e_2)(e_1-e_3)}{(a-e_1)^2}  \notag \\
& & 
\hspace{10mm}
+
\frac{t^2}{18}
+
\frac{t^2}{6}
\frac{e_1-e_3}{a-e_1}
+
\frac{t^2}{6}
\frac{e_1-e_2}{a-e_1}.
\end{eqnarray}
\par
From (\ref{eqn:A}), (\ref{eqn:B}),
(\ref{eqn:C}), (\ref{eqn:D}), 
we get
\begin{eqnarray*}
\Res_{x=e_1} 
\frac12 
\trace 
\left(
Y^{\prime}(x)Y^{-1}(x) 
\right)^2
&=& 
-
\frac{1}{24} 
\sum_{j \neq 1} 
\frac{1}{e_1-e_j}  \\
&+&
\frac
{
\sigma[p,q]^{\prime \prime}(t)
}
{\sigma[p,q](t)}  
\frac{1}{4(e_1-e_2)(e_1-e_3)} \\
&+&
t 
\frac
{
\sigma[p,q]^{\prime}(t)
}
{\sigma[p,q](t)}
\left( 
\frac16 \sum_{j \neq 1} 
\frac{1}{e_1-e_j} 
\right)
+ 
t 
\frac
{
\sigma[p,q]^{\prime}(t)
}
{\sigma[p,q](t)}
\frac12 \frac{1}{(a-e_1)}  \\
&+& 
\frac{t^2}{36}\frac{e_1-e_3}{e_1-e_2} 
+
\frac{t^2}{36}\frac{e_1-e_2}{e_1-e_3}  \\
&+& 
\frac{t^2}{4}
\frac{(e_1-e_2)(e_1-e_3)}{(a-e_1)^2} \\
&+& 
\frac{t^2}{18}
+
\frac{t^2}{6}\frac{e_1-e_3}{a-e_1}
+
\frac{t^2}{6}\frac{e_1-e_2}{a-e_1}.
\end{eqnarray*}
\par
By directly calculating the third term of 
$
\Res_{x=e_1} 
\frac12 
\trace 
\left(
Y^{\prime}(x)Y^{-1}(x) 
\right)^2,
$ 
we get
\begin{equation*}
t
\frac
{
\sigma[p,q]^{\prime}(t)
}
{
\sigma[p,q](t)
}
\left(
\frac16
\sum_{j \neq 1}
\frac{1}{e_1-e_j}
\right)
=
\frac{\eta_1 t^2}{6 \omega_1}
\left(
\frac{1}{e_1-e_2}
+
\frac{1}{e_1-e_3}
\right) 
+
\frac{t}{6 \omega_1}
\frac
{
\theta[p,q]^{\prime}
(\frac{t}{\omega_1}; \Omega)
}
{
\theta[p,q]
(
\frac{t}{\omega_1}; \Omega
)
}
\left(
\frac{1}{e_1-e_2}
+
\frac{1}{e_1-e_3}
\right).
\end{equation*}
From Lemma \ref{lem:sig}, 
(\ref{eqn:heat})
and (\ref{eqn:log}), 
we obtain 
\begin{eqnarray*}
\frac{\eta_1 t^2}{6 \omega_1}
\left(
\frac{1}{e_1-e_2}
+
\frac{1}{e_1-e_3}  
\right)  
&=&
\frac{t^2}{6}
\left(
-
\frac{1}{3 \omega_1^2}
\frac
{
\theta_{11}^{\prime \prime \prime}
}
{\theta_{11}^{\prime}}
\right)
\left(
\frac{1}{e_1-e_2}
+
\frac{1}{e_1-e_3}
\right)   \\
&=&
-
\frac{t^2}{3}
\frac
{
\frac{\partial \omega_1}{\partial e_1}
}
{\omega_1}
(e_1-e_2)
-
\frac{t^2}{3}
\frac
{
\frac{\partial \omega_1}{\partial e_1}
}
{\omega_1}
(e_1-e_3)
-
\frac{t^2}{9}
-
\frac{t^2}{18}
\frac{e_1-e_2}{e_1-e_3}
-
\frac{t^2}{18}
\frac{e_1-e_3}{e_1-e_2}.
\end{eqnarray*}
\par
By directly calculating the second term of 
$
\Res_{x=e_1} 
\frac12 
\trace 
\left(
Y^{\prime}(x)Y^{-1}(x) 
\right)^2,
$ 
\begin{eqnarray*}
& &
\frac
{
\sigma[p,q]^{\prime \prime}(t)
}
{\sigma[p,q](t)}
\frac{1}{4(e_1-e_2)(e_1-e_3)}
=
\left(
\frac{\eta_1}{\omega_1}
\right)
\frac{1}{4(e_1-e_2)(e_1-e_3)}
+
\left(
\frac{\eta_1 t}{\omega_1}
\right)^2
\frac{1}{4(e_1-e_2)(e_1-e_3)}  \\
& & \hspace{27mm} +
\frac{2 \eta_1 t}{\omega_1^2}
\frac
{
\theta[p,q]^{\prime}
(
\frac{t}{\omega_1}; \Omega
)
}
{
\theta[p,q](\frac{t}{\omega_1}; \Omega)
}
\frac{1}{4(e_1-e_2)(e_1-e_3)}
+
\frac{1}{\omega_1^2}
\frac
{
\theta[p,q]^{\prime \prime}
(
\frac{t}{\omega_1}; \Omega
)
}
{
\theta[p,q](\frac{t}{\omega_1}; \Omega)
}
\frac{1}{4(e_1-e_2)(e_1-e_3)}.
\end{eqnarray*}
From Lemma \ref{lem:sig}, 
(\ref{eqn:heat})
and (\ref{eqn:log}), 
each of four terms of 
$$
\frac
{
\sigma[p,q]^{\prime \prime}(t)
}
{\sigma[p,q](t)}
\frac{1}{4(e_1-e_2)(e_1-e_3)}
$$
is as follows:
\begin{eqnarray}
\label{eqn:A1}
\left(
\frac{\eta_1}{\omega_1}
\right)
\frac{1}{4(e_1-e_2)(e_1-e_3)}
&=&
-\frac12
\frac
{
\frac{\partial \omega_1}{\partial e_1}
}
{\omega_1}
-
\frac{1}{12}
\frac{1}{(e_1-e_2)}
-
\frac{1}{12}
\frac{1}{(e_1-e_3)},
\end{eqnarray}
and
\begin{eqnarray}
\label{eqn:B1}
& &
\left(
\frac{\eta_1 t}{\omega_1}
\right)^2
\frac{1}{4(e_1-e_2)(e_1-e_3)} \notag \\
&=&
t^2
\left(
\frac
{
\frac{\partial \omega_1}{\partial e_1}
}
{\omega_1}
\right)^2
(e_1-e_2)(e_1-e_3)
+
\frac{t^2}{36}
\frac{e_1-e_3}{e_1-e_2}
+
\frac{t^2}{36}
\frac{e_1-e_2}{e_1-e_3}  \notag \\
& &
+
\frac{t^2}{3}
\left(
\frac
{
\frac{\partial \omega_1}{\partial e_1}
}
{\omega_1}
\right)
(e_1-e_2)
+
\frac{t^2}{3}
\left(
\frac
{
\frac{\partial \omega_1}{\partial e_1}
}
{\omega_1}
\right)
(e_1-e_3)
+
\frac{t^2}{18},  
\end{eqnarray}
and
\begin{eqnarray}
\label{eqn:C1}
& &
\frac{2 \eta_1 t^2}{\omega_1^2}
\frac
{
\theta [p,q]^{\prime}
( \frac{t}{\omega_1} ; \Omega)
}
{
\theta [p,q]
(
\frac{t}{\omega_1}; \Omega
)
}
\frac{1}{4(e_1-e_2)(e_1-e_3)}  \notag \\
&=&
-
\frac
{
\frac{\partial \omega_1}{\partial e_1}
}
{\omega_1^2}
t
\frac
{
\theta[p,q]^{\prime}
(\frac{t}{\omega_1}; \Omega)
}
{
\theta[p,q](\frac{t}{\omega_1};\Omega)
}  \notag \\
& & -
\frac{t}{6 \omega_1}
\frac{1}{(e_1-e_2)}
\frac
{
\theta[p,q]^{\prime}(\frac{t}{\omega_1}; \Omega)
}
{
\theta[p,q](\frac{t}{\omega_1};\Omega)
}
-
\frac{t}{6 \omega_1}
\frac{1}{(e_1-e_3)}
\frac
{
\theta[p,q]^{\prime}
(\frac{t}{\omega_1}; \Omega)
}
{
\theta[p,q](\frac{t}{\omega_1};\Omega)
},
\end{eqnarray}
and
\begin{eqnarray}
\label{eqn:D1}
\frac{1}{\omega_1^2}
\frac
{
\theta[p,q]^{\prime \prime}
(\frac{t}{\omega_1};\Omega)
}
{
\theta[p,q](\frac{t}{\omega_1};\Omega)
}
\frac{1}{4(e_1-e_2)(e_1-e_3)}  
&=&
\frac
{
\frac{\partial}
{\partial \Omega} 
\theta[p,q](\frac{t}{\omega_1};\Omega)
}
{
\theta[p,q](\frac{t}{\omega_1};\Omega)
}  
\frac{\partial \Omega}{\partial e_1}.
\end{eqnarray}
We note that the sum of 
(\ref{eqn:C1}) 
and 
(\ref{eqn:D1})
is
\begin{equation*}
\frac{\partial}{\partial e_1}
\Big\{
\log
\theta[p,q]
\left(
\frac{t}{\omega_1}; \Omega
\right)
\Big\}
-
\frac{t}{6 \omega_1}
\frac{1}{(e_1-e_2)}
\frac
{
\theta[p,q]^{\prime}
(
\frac{t}{\omega_1};\Omega
)
}
{\theta[p,q](\frac{t}{\omega_1};\Omega)}
-
\frac{t}{6 \omega_1}
\frac{1}{(e_1-e_3)}
\frac
{
\theta[p,q]^{\prime}
(
\frac{t}{\omega_1};\Omega
)
}
{\theta[p,q](\frac{t}{\omega_1};\Omega)}.
\end{equation*}
\par
Therefore, 
we obtain
\begin{eqnarray*}
& &
\Res_{x=e_1}
\frac12
\left(
Y^{\prime}(x)
Y^{-1}(x)
\right)^2  \\
&=&
-
\frac18
\left(
\sum_{\mu \neq 1}
\frac{1}{e_1-e_{\mu}}
\right)
-
\frac12
\frac{\partial}{\partial e_1}
\log \omega_1  \\
&+&
t^2
\left(
\frac{\partial}{\partial e_1}
\log \omega_1
\right)^2
\prod_{\mu \neq 1}
(e_1-e_{\mu})
+
\frac{\partial}{\partial e_1}
\big\{
\log
\theta[p,q]
\left(
\frac{t}{\omega_1};\Omega
\right)
\big\}  \\
&+&
\frac{t}{2(a-e_1)}
\frac
{
\sigma[p,q]^{\prime}(\frac{t}{\omega_1};\Omega)
}
{
\sigma[p,q](\frac{t}{\omega_1};\Omega)
}
+
\frac{t^2}{4}
\frac
{
\prod_{\mu \neq 1}
(e_1-e_{\mu})
}
{(a-e_1)^2}
+
\frac{t^2}{6(a-e_1)}
\left(
\sum_{\mu \neq 1}
(e_1-e_{\mu})
\right).
\end{eqnarray*}

\end{proof}
We 
obtain $\omega_a$ 
in Proposition \ref{prop:omegaa} 
and 
get $\omega_{e_{\nu}}$ 
in Lemma \ref{lem:omeganu}. 
Therefore, we have $H_{\nu} \, (\nu=1,2,3)$ 
in the following way:
\begin{proposition}
For $\nu=1, 2, 3,$
\begin{eqnarray*}
H_{\nu}
&=&
\frac{\partial }{\partial e_{\nu}}
\Big\{ 
\log 
\theta [p,q] 
\left( 
\frac{t}{\omega_1} ; \Omega 
\right) 
\Big\}
- 
\frac12 
\frac{\partial }{\partial e_{\nu}}
 \big\{ \log \omega_1 \big\}
-
\frac18 
\sum_{\mu \neq \nu} \frac{1}{(e_{\nu}-e_{\mu})}  \\
&+& 
\frac{\partial}{\partial e_{\nu}}
\left(
\frac{\eta_1 t^2}{2 \omega_1}
\right)
+
\frac{\partial}{\partial e_{\nu}} 
\Big\{ \frac{t^2}{4} f(e_1, e_2, e_3) \Big\}.
\end{eqnarray*}
\end{proposition}

Since we have calculated the Hamiltonians 
$H_t, H_1, H_2, H_3$, 
we finally get the following theorem:
\begin{theorem}
For the monodromy preserving deformation (\ref{eqn:mpd}),
the $\tau$ function is
\begin{equation*}
\tau(e_1, e_2, e_3 ; t)
= 
\theta[p,q]
\left(
\frac{t}{\omega_1};\Omega
\right)
\omega_1^{-\frac12}
\prod_{1 \leq \nu<\mu \leq 3} 
(e_{\nu}-e_{\mu})^{-\frac18}   
\exp(\frac{\eta_1 t^2}{2 \omega_1})
\exp(\frac{t^2}{4} f(e_1, e_2, e_3)),
\end{equation*}
where 
\begin{equation*}
f(e_1, e_2, e_3) 
= 
-a 
+ 
\frac13 \sum_{\nu=1}^3 e_{\nu}  
+ 
\frac12 \prod_{\nu=1}^3(a-e_{\nu})
\sum_{\mu=1}^3 \frac{1}{(a-e_{\mu})^2}.
\end{equation*}
\end{theorem}

\appendix
\section{Appendix}
In \cite{MM0}, Miwa-Jimbo expressed 
the $\tau$-function 
with the Weierstrass $\sigma$-function. 
The paper, however, is written in Japanese. 
We will write the summary in English. 
Here, we assume 
$e_1+e_2+e_3=0$
\par
We fix a point $\alpha$ in one dimension complex 
torus $\T^1$ with period 
$\omega_1, \omega_2.$ 
We take $l \in \Z$ 
and 
the parameter $t \in \C.$ 
We consider the following function:
\begin{eqnarray*}
y_l(z) 
&=&
y_l(z,t,\alpha) \\
&=&
\sigma(z-\alpha)^{l-1}
\sigma(z+\alpha)^{-l}
\sigma(z+t+(2l-1)\alpha)
\exp 
\{
-
\frac{t}{2}
\left(
\zeta(u-\alpha)
+
\zeta(u+\alpha)
\right)
\}
\end{eqnarray*}
From the quasi-periodicity of the $\sigma$-function 
and the $\zeta$-function, 
$y_l(z)$ is holomorphic in $\T$ 
except for $\alpha$. Near $z=\alpha$, 
\begin{eqnarray*}
& &
y_l(z)
=
\hat{y}_l(z)
(z-\alpha)^{l-1}
\exp
\left(
-
\frac{t}{2}
\left(
\frac{1}{z-\alpha}
-
\frac
{
\wp^{\prime \prime}(\alpha)
}
{
2
\wp^{\prime}(\alpha)
}
\right)
\right)\\
& &
\hat{y}_l(z)
=
c_0
\left(
1
+
c_1(z-\alpha)
+
\cdots
\right) \\
& &
c_0
=
\sigma(2\alpha)^{-l}
\sigma(t+ 2l \alpha)
\exp 
\left(
-
\frac{t}{2}
\zeta(2 \alpha)
-
\frac{t}{4}
\frac
{
\wp^{\prime \prime}(\alpha)
}
{
\wp^{\prime}(\alpha)
}
\right)
=
c_0(t, \alpha, l)\\
& &
c_1
=
\zeta(t+2l \alpha)
-
l
\zeta(2\alpha)
+
\frac{t}{2}
\wp(2\alpha)
=
c_1(t, \alpha, l).
\end{eqnarray*}
The above property characterizes $y_l(z)$.
\par
We set 
$e_1, e_2, e_3, e_{\infty}=\infty \in \bbP^{1}$ 
as 
the branch points of
the covering map 
$
\T \rightarrow \bbP^{1}, 
z \rightarrow x=\wp(z).
$
We define the multi-valued analytic matrix $Y(x)$ 
in the following way:
\begin{eqnarray*}
& &
Y(x)
=
G_0^{-1}
\tilde{Y}(z) \\
& &
G_0^{-1}
=
\diag
\left(
\frac{1}{c_0(t,\alpha,l)},
\frac{1}{c_0(-t,\alpha, -l)}
\right)\\
& &
\tilde{Y}(z)
=
\left(
\begin{array}{cc}
y_l(z) & y_l(-z) \\
y_{l+1}(z) & y_{l+1}(-z) 
\end{array}
\right).
\end{eqnarray*}
We assume $\wp(\alpha)\neq e_1, e_2, e_3, e_{\infty}.$
Then 
$Y(x)$ 
satisfies 
the following linear ordinary differential equation:
\begin{eqnarray*}
& &
\frac{dY}{dx}=A(x)Y \\
& &
A(x)
=
\frac{A_{-2}}{(x-a)^2}
+
\frac{A_{-1}}{x-a}
+
\frac{B_1}{x-e_1}
+
\frac{B_2}{x-e_2}
+
\frac{B_3}{x-e_3}\\
& &
A_{-2}
=
\left(
\begin{array}{cc}
\frac{\bar{t}}{2} &                   \\
                  & -\frac{\bar{t}}{2}
\end{array}
\right),
\bar{t}=\wp^{\prime}(\alpha)t,
a=\wp(\alpha).
\end{eqnarray*}
In order to check this, 
we will study the monodromy of $Y(x)$.

\subsection{The Behavior near the Irregular Singular Point}
We note that 
\begin{equation}
z
-
\alpha
=
\frac{1}{\wp^{\prime}(\alpha)}
(x-a)
-
\frac{\wp^{\prime \prime}(\alpha)}
{
2 \wp^{\prime \prime \prime}(\alpha)
}
(x-a)^2
+
\left(
\frac{\wp^{\prime \prime}(\alpha)^2}
{
2 \wp^{\prime}(\alpha)^5
}
-
\frac
{
\wp^{\prime \prime \prime}(\alpha)
}
{6 \wp^{\prime}(\alpha)^4}
\right)
(x-a)^3
+
\cdots,
\end{equation}
We get the following equations:
\begin{eqnarray*}
& &
Y(x)
=
\hat{Y}^{(\alpha)}(x)
\exp 
\{
T^{(\alpha)}(x)+K
\} \\
& &
T^{(\alpha)}(x)
=
\left(
\begin{array}{cc}
\frac{\bar{t}}{2} &                   \\
                  & -\frac{\bar{t}}{2}
\end{array}
\right)
\frac{-1}{x-a}
+
\left(
\begin{array}{cc}
l-1 &   \\
    &-l-1 
\end{array}
\right)
\log(x-a) \\
& &
K
=
\left(
\begin{array}{cc}
-l+1 & \\
     & l+1 
\end{array}
\right)
\log \wp^{\prime}(\alpha) \\
& &
\hat{Y}^{(\alpha)}(x)
=
\left(
1
+
Y^{(\alpha)}_1(x-a)
+
\cdots
\right) \\
& &
Y_1^{(\alpha)}
=
\frac{1}{\wp^{\prime}(\alpha)}
\tilde{Y}_1^{(\alpha)}
+
\frac{t}{2}
\left(
\frac
{
\wp^{\prime \prime}(\alpha)^2
}
{
4\wp^{\prime}(\alpha)^3
}
-
\frac
{
\wp^{\prime \prime \prime}(\alpha)
}
{
6 \wp^{\prime}(\alpha)^2
}
\right)
\left(
\begin{array}{cc}
1 & \\
  & -1
\end{array}
\right)
-
\frac
{
\wp^{\prime \prime}(\alpha)
}
{
\wp^{\prime}(\alpha)^2
}
\left(
\begin{array}{cc}
l-1 & \\
    & -l-1
\end{array}
\right) \\
& &
\hat{Y}_1^{(\alpha)}
=
\left(
\begin{array}{cc}
c_1(t,\alpha,l) 
& 
\frac{c_0(-t,\alpha,1-l)}{c_0(t,\alpha,l)} \\
\frac{c_0(t,\alpha,l+1)}{c_0(-t,\alpha,-l)} 
& 
c_1(-t,\alpha,-l)
\end{array}
\right).
\end{eqnarray*}

\subsection{The Behavior near Regular Singular Points}
We assume $\nu=1,2,3.$ By direct calculation, we get 
\begin{eqnarray*}
& &
y_l(z;t,\alpha)
=
y_l(\tilde{\omega}_{\nu})
(1+f_l^{(\nu)}(z-\tilde{\omega}_{\nu})+\cdots) \\
& &
y_l(\tilde{\omega}_{\nu};t,\alpha)
=
\frac
{
\sigma_{\nu}(t+(2l-1)\alpha)
}
{\sigma_{\nu}(\alpha)}\\
& &
f_l(\tilde{\omega}_{\nu})
=
\zeta_{\nu}(t+(2l-1)\alpha)
-
\zeta_{\nu}(\tilde{\omega}_{\nu}+\alpha)
+
t\wp(\alpha+\tilde{\omega}_{\nu}),
\end{eqnarray*}
where 
\begin{eqnarray*}
& &
\sigma_{\nu}(z)
=
\exp 
\{
-
\frac{\tilde{\eta_{\nu}}}{2}
z
\}
\frac
{
\sigma(z+\tilde{\omega}_{\nu})
}
{
\sigma(\tilde{\omega}_{\nu})
},\\
& &
\zeta_{\nu}(z)
=
\zeta(z+\tilde{\omega}_{\nu})
-
\frac
{
\tilde{\eta_{\nu}}
}
{2},
\end{eqnarray*}
and $\tilde{\omega}_{\nu}, \tilde{\eta_{\nu}}$ 
are defined in section 5.
We set $G^{(\nu)}$ in the following way:
\begin{equation*}
Y(x)
=
G^{(\nu)}
\hat{Y}^{(\nu)}(x)
(x-e_{\nu})^
{
\left(
\begin{array}{cc}
0 & \\
  & \frac12
\end{array}
\right)
}
\left(
\begin{array}{cc}
1 & 1 \\
1 & -1
\end{array}
\right).
\end{equation*}
And
\begin{eqnarray*}
G^{(\nu)}
\hat{Y}^{(\nu)}(x)
&=&
G_0^{-1}
\frac12
\left(
\begin{array}{cc}
y_l(z)
+
y_l(-z) 
& 
\frac
{y_l(z)-y_l(-z)}
{z-\tilde{\omega}_{\nu}} \\
y_{l+1}(z)
+
y_{l+1}(-z) & 
\frac
{y_{l+1}(z)-y_{l+1}(-z)}
{z-\tilde{\omega}_{\nu}} 
\end{array}
\right)
\left(
\begin{array}{cc}
1 & \\
  & 
\frac{z-\tilde{\omega}_{\nu}}{\sqrt{x-e_{\nu}}}
\end{array}
\right) \\
&=&
G^{(\nu)}
(1+Y_1^{(\nu)}(x-e_{\nu})+\cdots) \\
G^{(\nu)}
&=&
G_0^{-1}
\left(
\begin{array}{cc}
y_l(\tilde{\omega}_{\nu}) & \\
               & y_{l+1}(\tilde{\omega}_{\nu})
\end{array}
\right)
\left(
\begin{array}{cc}
1 & f_l^{(\nu)} \\
1 & f_{l+1}^{(\nu)}
\end{array}
\right)
\left(
\begin{array}{cc}
1 & \\
  & 
\sqrt
{
\frac{2}
{
\wp^{\prime \prime}(\tilde{\omega}_{\nu})
}
}
\end{array}
\right)
\end{eqnarray*}
\par
The case $x=\infty, z=0$ can be calculated 
in the same way. 
By direct calculation, we get 
\begin{eqnarray*}
& &
y_l(0;t,\alpha)
=
(-1)^{l-1}
\frac
{
\sigma(t+(2l-1)\alpha)
}
{
\sigma(\alpha)
} \\
& &
f_l^{(\infty)}
=
\zeta(t+(2l-1)\alpha)
-
(2l-1)\zeta(\alpha)
+
t\wp(\alpha).
\end{eqnarray*}
And
\begin{eqnarray*}
& &
Y(x)
=
G^{(\infty)}
\hat{Y}^{(\infty)}(x)
\left(
\frac{1}{x}
\right)^
{
\left(
\begin{array}{cc}
0 & \\
  & \frac12
\end{array}
\right)
}
\left(
\begin{array}{cc}
1 & 1 \\
1 & -1
\end{array}
\right)\\
& &
\hat{Y}^{(\infty)}(x)
=
1+O
\left(
\frac{1}{x}
\right) \\
& &
G^{(\infty)}
=
G_0^{-1}
\left(
\begin{array}{cc}
y_l(0) & \\
       & y_{l+1}(0)
\end{array}
\right)
\left(
\begin{array}{cc}
1 & f_l^{(\infty)}  \\
1 & f_{l+1}^{(\infty)}
\end{array}
\right).
\end{eqnarray*}

\subsection{The tau-function}
According to the definition of 
the $\tau$-function, 
\begin{eqnarray*}
\frac{\partial}{\partial t}
\tau(t)
&=&
\wp^{\prime}(\alpha)
\trace 
Y_1^{(\alpha)}
\left(
\begin{array}{cc}
\frac12 &   \\
        & -\frac12
\end{array}
\right)\\
&=&
\zeta(t+2l \alpha)
+
\frac{t}{2}
\left(
\wp(2\alpha)
+
\frac
{
\wp^{\prime \prime}(\alpha)^2
}
{
4 \wp^{\prime}(\alpha)^2
}
-
\frac
{
\wp^{\prime \prime \prime}(\alpha)
}
{
6\wp^{\prime}(\alpha)
}
\right)\\
& &
\hspace{20mm}
-l
\left(
\zeta(2\alpha)
+
\frac
{
\wp^{\prime \prime}(\alpha)
}
{
2\wp^{\prime}(\alpha)
}
\right).
\end{eqnarray*}
And
\begin{eqnarray*}
& &
\tau(t)
=
\tau_{l}(t)
=
\sigma(t+2l\alpha)
\exp h_l(t) \\
& &
h_l(t)
=
\frac{t^2}{4}
\left(
\wp(2 \alpha)
+
\frac
{
\wp^{\prime \prime}(\alpha)^2
}
{
4\wp^{\prime}(\alpha)^2
}
-
\frac
{
\wp^{\prime \prime \prime}(\alpha)
}
{
6\wp^{\prime}(\alpha)
}
\right)
-
tl
\left(
\zeta(2\alpha)
+
\frac
{
\wp^{\prime \prime}(\alpha)
}
{
2\wp^{\prime}(\alpha)
}
\right).
\end{eqnarray*}

\end{document}